\documentstyle[12pt]{article}
\begin{document}

\title{Complex  $su_{q}(2)$ dynamical symmetry, limiting
cases and  one-dimensional 
potential realisations}

\author{A. Ludu$^1$ and W. Scheid  \\
\normalsize{Institut f\"ur Theoretische Physik der } \\ 
\normalsize{Justus-Liebig-Universit\"at, 
Heinrich-Buff-Ring 16, D-35392, Giessen, Germany}
}
\date{}
\maketitle
\begin{abstract}
Using a complex deformation ($q=e^{is}$) of the Lie algebra $su(2)$
we obtain extensions of the finite-dimensional
 representations towards
infinite-dimensional ones.
A generalised q-deformation of $su(2)$, as a
Hopf algebra, is
introduced and 
it is proved that some infinite-dimensional representations
 can be also constructed.
The Schr\"{o}dinger
picture of
$su_{e^{is}}(2)$ is investigated, 
using a differential realization, and a large class of 
equivalent potentials is obtained. Limiting cases
of the deformations  are also analysed.
A connection between the case when $q$ is a root of $1$ and the
commensurability of the corresponding potentials is found.

\end{abstract}
\setlength{\baselineskip} {4ex}.

\footnotetext[1]{Permanent Adress: 
Bucharest University, Department of Theoretical Physics, Bucharest-M\u agurele
P.O. Box MG-5211, Romania,
e-mail: ludu@th.physik.uni-frankfurt.de} 
\vskip2cm

{\bf Short title:} Complex q-deformation of su(2) and  potential realization
\vskip 2truecm
{\bf PACS:} 03.65.Fd, 11.30.Na 
\vfill
\eject

\section{Introduction}

\vskip 1truecm
Quantized universal enveloping algebras (also called q-algebras, q-deformed
Lie algebras)
 have been the subject of numerous
recent studies in mathematics and
physics. They represent some special deformations ($q\neq 1$) of the 
universal enveloping algebra of Lie algebras ($q=1$) [1].
They were first introduced  
by Sklyanin [2] and by Kulish and Reshetikhin [3],
and
developed through a Jordan-Schwinger realisation of $su_{q}(2)$ by 
Biedenharn [4], Macfarlane [5] and Ng [6], or through $su_{q}(1,1)$ by
Ui and Aizawa [7], Kulish and Damaskinsky [8], Maekawa 
and Masuta {\it et al} [9] and Gromov and
Man'ko [10].
Since then $su_{q}(2)$ and $su_{q}(1,1)$ have been applied in various
branches of physics.
The importance of these deformations, in one-dimensional problems,
establishes its role beyond doubt for many physical problems, like:
deformed spin-chain models [11], description of rotation-vibration
spectra of deformed or superdeformed nuclei [12-14], molecular spectra [15],
as well as two-dimensional $su_{q}(3)$ [16] and three-dimensional $u_{q}(6)$
[17], or $su_{q}(N)$ [18] generalisations.
\par
Such deformations provide much interest in various other
contexts, e.g. q-deformed
SUSY [19], exact solvable potentials [20], Hamiltonian quantization [21]
and fractional anyonic statistics [22]. Consequently they might prove useful in
applications to nuclear, molecular physics [23] and scattering theory.
The method of equivalent potentials developed in [24] sugests the possibility
of a continuous (q-deformation) of some dynamical symmetries,
i.e. bound states towards scattering states by deforming
the harmonic
oscillator potential into a P\"oschl-Teller-like potential, by both
real and complex 
q-deformation
of the energy levels.
We also note that q-deformation {\bf with complex q} have been investigated
and proved to be usefull in many physical applications [9,10,13,17,22].
\par
Recently, new generalised forms of deformation of Lie algebras have been 
introduced: deformations involving $J_{0}$ only in one commutator
[25,26] (the latter reference proving a generalised deformation endowed with a
Hopf algebraic structure), or in all theree commutators [27].

A possible procedure, inspired by the generalization
of the Algebraic Scattering Theory (AST)  into q-deformations, could be:
to q-deform a given symmetry with physical applications, 
and to smoothly modify the q parameter,
in such a way, to find, for certain limiting values, new symmetries
of physical interest. More than this, one can investigate what is 
happening in between
these limiting cases, and the possible underlying physics (
and some physical significance of the parameter $q$).
\par
In this paper we add to this field and prove that some certain
q-deformations
(trigonometric in Sect.2 or more general, in Sect.3) 
can carry the bounded unirreps of $su(2)$ to
unbounded unitary representations, some of them being equivalent to those of $su(1,1)$.
We also prove the equivalence of $su_{q}(2)$ with the two-dimensional
Euclidean Lie algebra $e(2)$ for certain values of $q$.
In Sect.4 we find a differential realisation of $su_{q}(2)$ which can be 
related to a
Schr\"odinger equation coresponding to some particular one-dimensional
potentials of physical interest. 
These results prove that structures more general then Lie
algebras could succesfully
be applied to algebraic scattering theory (AST) [28]. 
\par
The aim of this paper is to implement the q-deformation in
the traditional method used in
AST, i.e. starting from systems
with bound state solutions, described by a compact
algebra, and continuosly deforming in order to obtain scattering solutions.
In traditional AST this is achived by analytic continuation
into the complex plane of the generators (the compact Lie algebra is 
analytically continued into a
non compact one by imposing supplementary conditions on the generators).
In the present paper we prove that this 
could be obtained
by using the q-deformation of the compact $su(2)$ Lie algebra
to $su_{q}(2)$, $q=e^{is}$, the later having infinite dimensional 
unitary representations.
In this way the transition from bounded to scattering states 
in one-dimensional systems could be
described in a unified algebraic formalism.
More than these, a possible way to get the Euclidean connection [28]
through q-deformation, by starting from the compact Lie algebra $su(2)$
and obtain a q-deformed algebra similar with $e(2)$ is worked out in Sect.2.
This algebraic result was also realised in a potential picture with
q-differential operators for the generators.
A large set of potential shapes were obtained from the equivalent
Schr\"odinger equation, for different values of $q$. We have found a
similarity between the rational/irrational character of $s/ \pi $ and the 
random/incommensurate character of certain exact integrable potentials.

The other possibility, we have found, is to work with $q\in R$ and to modify
the defining function $[2J_{z}]$. Such an example is worked out in Sect.3,
where the non-linear algebra has also a Hopf algebra structure
[1,29-32]. This could
be usefull in further applications in the quantum mechanics algebraic
approach of many
particles systems and lattice models.

\vskip 1.5truecm
\section{Unitary representations of $su_{q=e^{is}}(2)$}
\vskip 1truecm
\par

The purpose of this section is to give a brief discussion of some
complex q-deformations   
of the Lie algebra $su(2)$ together with some of its unitary representations.
In the following we review some results concerning the                      
q-algebra $su_{q}(2)$, $q$ generic, [4-5,9,10,29,30].
In the $su_{q}(2)$, generated by the 
operators 
$J_{+}$, $J_{-}$ and $J_{z}$, one has the following commutation relations
\begin{eqnarray}
[J_{z},J_{\pm }] & =& \pm J_{\pm }, \hfill
\end{eqnarray}
\begin{eqnarray}
[J_{+},J_{-}] & =& {
{ q^{2J_{z}}-q^{-2J_{z}}
} 
\over {q-q^{-1}}
}
=[2J_{z}]_{q}=[2J_{z}], \hfill
\end{eqnarray}
where $J_{\pm}=J_{x} \pm iJ_{y}$ ($J_{x}$, $J_{y}$ and $J_{z}$ self-adjoint
operators). 
The structure given
by eqs.(1,2) can be endowed with a real Hopf algebra structure 
[1,29-30].
\par
The Casimir operator C is given by
\begin{eqnarray}
C & = & \biggl [ J_{z}\pm {1 \over 2}\biggr ] ^{2}+
J_{\mp}J_{\pm}={{[2]} \over 2}[ J_{z}]^{2}+{1 \over 2}
(J_{+}J_{-}+J_{-}J_{+})+\biggl [{1 \over 2}\biggr ]^{2}. \hfill
\end{eqnarray}
For $q \in R$ the  unirreps of $su_{q}(2)$, all finite-dimensional,
acting on a Hilbert space $V_{j}$ with the basis $|jm>$, are 
characterized by the relations (e.g. [5])
\begin{eqnarray}
C|jm> & = & \biggl [j+{1 \over 2}\biggr ]^{2}|jm>, \hfill
\end{eqnarray}
\begin{eqnarray}
J_{z}|jm> & = & m|jm>,\hfill
\end{eqnarray}
\begin{eqnarray}
J_{\pm}|jm> & = & ([j \mp m][j \pm m+1])^{1/2}|jm \pm 1>,\hfill
\end{eqnarray}
where m takes integers or half-integers ranging from $-j$ 
to $j$ for each
$j=0,{1 \over 2},1,...$. These representations are in one-to-one relation
with the unirreps of the Lie algebra $su(2)$ and coincide with them in the
limit $q \rightarrow 1$.
Infinite dimensional representations for the quantum
algebra $su_{q}(1,1)$
can be obtained  in a similar way. 
The commutators in eq.(1) remain unchanged and eq.(2) becomes
\begin{eqnarray}
[J_{+},J_{-}] & =& -[2J_{z}].  \hfill
\end{eqnarray}
One obtains the four classes of unirreps for $su_{q}(1,1)$ 
 [9,11, 21, 29-31] 
(the principal 
and complementary
series for continuous variation of $j$ and positive and negative discrete series
for discrete $j$) labeled by $j$, through
\begin{eqnarray}
C|jm> & =& [j][j \pm 1]|jm>.  \hfill
\end{eqnarray}
The unirreps reduce also to those of $su(1,1)$,
 in the limit $q \rightarrow 1$.
We mention that several excellent papers about this topics exist (i.e.
unirreps of $su_{q}(2)$ and
$su_{q}(1,1)$, or, in general
$su_{q}(n)$ or
$su_{q}(n,m)$,  [8-10,21,29-31], and for a recent review see [30] and
references
herein,
which analyse and clasify their unirreps. Here, we only want
to give a short review and  pictural examples which can exemplify the status
of the problem and could be used to give an easier way of understanding and
using these techniques.

It is already a classical result [1,5-7,21,29,30-32] that, for real $q$, the
unirreps of $su_{q}(2)$ are equivalent with those of $su(2)$ (Lusztig and
Rosso's Theorem, see for details [29,30]).
This fact is remarkably explained by Biedenharn [29] as being the
consequence of the invariance of the (discrete) integer number which labels
the dimension of the unirrep, versus the continuous variation of $q$,
especially in the case of simple compact Lie groups. In the case when $q$ is
a root of $1$, the representation theory becomes strikingly different due
to the fact that the generators and the Casimir operator become nilpotent
elements. The representations are indecomposable and no longer irreductible
in general [29,30]. These representations can be clasified, acording to
their q-dimension $d_{q}=\sum_{rep.}q^{J_{z}}$ into two classes:
with $d_{q}=0$ and with $d_{q}\neq 0$.

In the following we show that infinite-dimensional unitary representations
can also be obtained directly  from $su_{q}(2)$ with $q=e^{is}$.
We extend the q-deformed
algebra, $su_{q}(2)$,
eqs.(1,2), to a complex
q-deformed algebra, by imposing $q=e^{is}$, $s \in (0,\pi )$ 
(i.e. $|q|=1$) with $s\neq k\pi$, $k\in {N}$.
We discuss the cases $s \rightarrow 0$, $\pi /2$ and $\pi $
separately.
This algebra is not homeomorphic with $su_{q}(1,1)$ due to eq.(2)
but some similarities with the unbounded representations of
$su_{q}(1,1)$ and $sl_{q}(2)$, [33], are existing.
The equivalence between $su_{-1}(2)$ with $su(1,1)$, already studied
in [34] and noted in the conclusions of [22] represents the basis for our
following calculations.

In this Section $q=e^{is}$ and
\begin{eqnarray}
[x]_{q=e^{is}} & =& {{q^{x}-q^{-x}} \over {q-q^{-1}}}={{sin(xs)} \over {
sin(s)}}\equiv [x], \hfill
\end{eqnarray}
for a quantity $x$ (a number or operator) will be used without any label,
for $q=e^{is}$.
When $x$ is a real number, $[x]$ is also real.
We note that for $q=e^{is}$ the sign of $[x]$ is not necesary the
same with the sign of $x$, like it is always in the real q-deformation
case, $q=e^{s}$. 
In addition, for any sequence $\{ x_{n} \} _{n\in Z}$, the q-deformed
one $\{ [x]\} _{n\in Z}$ is always bounded.
These are the points which directed our studies
to investigate the $su_{q}(2)\simeq su(1,1)$ similarity.
For example, for a complex number $x=\alpha +i\beta $ we have
\begin{eqnarray}
[x] & = &[\alpha ]\sqrt{1+[\beta]_{e^{s}}^{2}
\biggl (1-{{[2]^{2}_{e^{s}}} \over 4} \biggr )}
+[i][\beta ]_{e^{s}}\sqrt{1- [\alpha]^{2}
\biggl (1-{{[2]^{2}} \over 4}\biggr )}. \hfill
\end{eqnarray}
The object $[x]$ is invariant under the transformation $q=e^{is}
\rightarrow q^{-1}=q^{*}$ or, equivalentely, $s \rightarrow -s$.
 
For unitary
representations we need the conditions $J_{\pm }^{\dagger }=J_{\mp }$ and
$J_{z}^{\dagger }=J_{z}$ which imply that the eigenvalues of $J_{+}J_{-}$
and $J_{-}J_{+}$ must be positive and those of $J_{z}$ must be real.
By computing the matrix elements of $J_{+}$ and by using eq.(2) it results
that two consecutive values of m must differ only by 
$\pm 1$, like in the classic case of $su(2)$. 
The operators $J_{z}$ and $C$ (eq.(3)) form 
a complete set of commuting self-adjoint operators, $C|cm>=c|cm>$, $J_{z}|cm>
=m|cm>$, $m \in R$ and $c\ge 0$.

As $|cm>$ is eigenvector for $J_{\pm}J_{\mp}$ we have
\begin{eqnarray}
J_{\pm }J_{\mp }|cm> & =& \biggl (c-\biggl [m\mp {1 \over 2}
\biggr ]^{2} \biggr )|cm>. \hfill
\end{eqnarray}
From eq.(11), the unitarity condition is fulfiled if (see for example
[35]):
\begin{eqnarray}
c \cdot \sin ^{2}s & \ge & \sin ^{2}s\biggl (m\mp 
{1 \over 2}\biggr ). \hfill
\end{eqnarray}
We consider now the unitary representations of $su_{e^{is}}(2)$, 
with $s\neq 0, \pi$
to be unitary extensions of the unirreps of $sl_{q}(2)$.
In both these limits $s \rightarrow 0$, $\pi $, eq.(12) gives the
finite-dimensional unirreps of $su(2)$. For $s \rightarrow 0$, $m$ is integer
or half-integer but for the limiting case $s \rightarrow \pi$,
$m$ is only integer. In the following we want to present some of the unitary
representations of $su_{e^{is}}(2)$ according to the general procedure used
in AST references. By this we mean that we follow the criterium 
of unitarity, given by
eq.(12), regarding the algebra defined with eqs.(1,2) with $q=e^{is}$.
We take here into account the algebraic structure, ignoring the Hopf
algebraic structure. This will be analysed, for other deformations, in
Section 3. We classify the unitary representations
by the different admissible ranges of $c$.
We do not take into account other existing equivalent representations,
like , for example, the one generated by the automorphism of
$su_{e^{is}}(2)$:
$J_{+} \rightarrow -J_{+}$, $J_{z} \rightarrow J_{z}+{{i \pi} \over 2}$,
without a classical counterpart [1,9,10,30-32]. This representation
can be obtained in our case by a shift of $\pi $ in $s$.
 
We investigate such unitary representations, which fall
into one of the following three categories: 
\vskip 0.7truecm
\underline{1. Infinite dimensional continuous series, Fig.1a} 

(c continuous, bounded from below, m unbounded)
\vskip 0.7truecm
We can take $s \in (0,\pi )$ since, for the other half of the interval,
$(\pi , 2\pi )$, the calculations are similar. We label the unirreps
by the $c$ and $m$ numbers. 
Following the inequality in eq.(12), this class, \underline{(1)},
is defined for $c$ within the range:
${1 \over {sin^{2}(s)}} < c < \infty $ 
(having continuous variation bounded from below) 
and m taking all the integers or 
half-integers (unbounded form above and below) 
with the restrictions $m_{k}\neq {\pi \over {2s}}\pm
{1 \over 2}+{{k\pi } \over s}$, ($m=m_{k}$ only in the limit $s=\pi /2$). 
We can also label these unitary representations
by corresponding values for $j={{\pi }\over s}\biggl (k+{1 \over 2}
\biggr )-{1 \over 2}+i\sigma$, $\sigma \in R$, $k \in Z$, i.e.:
\begin{eqnarray}
c =\biggl [j+{1 \over 2}\biggr ]^{2} & = &
{{cosh^{2}s\sigma} \over {\sin ^{2}s}} >{1 \over {sin ^{2}s}}. \hfill
\end{eqnarray}
One can see that the expression of the label $j$ is not uniquely determined
by $\sigma $, hence we have obtained the same uncertainty regarding
the labeling of these representations, like in the case when $q$ is 
a root of $1$, [29-30].
In these series the minimum value of the eigenvalue of 
$C$ coincides with the extreme
values of the Casimir operator of the "boundary" representations of $sl_{q}(2)$
(q a root of $1$) obtained by
Alexeev {\it et al}, (eq.(18) in [33]).
These series are also similar with the 
principal series in [9]
for $su_{q}(1,1)$ or with those belonging to the quantum Cayley-Klein algebras 
$su_{q}(2;j_{1})$ in [10] for $j_{1}=i$ ("the strange series").
The action of $J_{\pm}$ gives 
$$
J_{\pm }|c,m>=N_{cm}^{\pm }|c,m>
$$
\begin{eqnarray}
 & = & {{(\cos ^{2}{\beta }_{\pm }\cosh^{2}
(\sigma s)+\sin ^{2}{\beta }_{\pm }\sinh ^{2}(\sigma s))^{1/2}}
\over {\sin s}}|c,m\pm 1>,  \hfill
\end{eqnarray}
with ${\beta }_{\pm}=\mp ms -{s \over 2}$. The coefficients $N_{cm}^{\pm }$
are real and become zero when $\sigma \rightarrow 0$ and 
$m=\pm \biggl (k+{1 \over 2}\biggr ){{\pi }
\over s}\pm {1 \over 2}$. When $\sigma \rightarrow 0$ we obtain the limiting 
value for the
Casimir $c\rightarrow {1 \over {\sin ^{2}s }}$. Consequently, we have shown
that the representations of this class are unitary and infinite-dimensional,
according to eqs.(12,14). Since $N^{\pm }_{cm}\neq 0$,
it results that $J_{\pm }$ do not act nilpotently on the space of
representation and this becomes irreducible.
An example of such series is ploted in Fig.1a, with full circles, in the
range $c>c_{0}$.

\vskip 0.7truecm
\underline{2. Mixed series,
infinite and finite-dimensional unirreps, Fig.1a}.

(c continuous, bounded from below and above, m bounded or unbounded)
\vskip 0.7truecm
The alowed ranges of $c$ and $m$ for these series are
obtained again from eq.(12).
Due to the trigonometric
functions involved in this inequality, the allowed values for $m$
are in certain periodically distributed intervals, which we describe below.
In this case the Casimir eigenvalues take values in the range
$c_{1}={1 \over {4\sin ^{2}(s/2)}}<c\leq c_{0}={1 \over {\sin ^{2}s}}$.
The allowed range for $m$ is given by
\begin{eqnarray}
{\bar J}_{\delta }\cup {\bar J}_{\Delta}&=&
\biggl [{{\pi } \over {2s}}+k{{\pi} \over {s}}-{{\delta}
\over {2}},{{\pi } \over {2s}}+k{{\pi} \over {s}}+{{\delta}
\over {2}} \biggr ]\cup \biggl [{{\pi } \over {s}}+k{{\pi} \over
{s}}-{{\Delta }
\over {2}},{{\pi } \over {s}}+k{{\pi} \over {s}}+{{\Delta }
\over {2}} \biggr ] \hfill
\end{eqnarray} 
with $k \in {\bf Z}$ and
$J_{\delta } $, $J_{\Delta }$ being open intervals on the $m$-axis,
of lenghts
 $\delta={{s-\pi+2\alpha } \over {s}}$,
$\Delta={{2\alpha -s} \over {s}}$, respectively,
 with $\alpha=\arcsin \sqrt{c\sin ^{2}s}$.
We have denoted with ${\bar J}_{\delta }, {\bar J}_{\Delta}$ the closure of these
intervals.
These series allow the $\alpha $, $\delta $ and $\Delta $
parameters to take values in the intervals:
$\alpha \in \biggl [ {{\pi -s} \over {2}},{{\pi } \over 2}\biggr ]$,
$\delta \in [0,1]$ and $ \Delta \in \biggl [ 
{{\pi -2s} \over {s}},{{\pi -s } \over s}\biggr ] $, function of
the eigenvalue $c$.
In between the allowed intervals for $m$, there are forbiden intervals $J_{f}$
of equal lenght ${{\pi -2\alpha } \over {s}}\in [0,1]$. All these 
sequences ($J_{\delta }, J_{f}, J_{\Delta}, J_{f}$)
have periodicity ${{\pi } \over s}$. These series  can be divided into three
different sub-classes:
\vfill
{\it a. Infinite dimensional unitary representations}

(c continuous and bounded, m unbounded)

We obtain infinite dimensional unitary representations with bounded values for $c$ and
unbounded for $m$ when
$s={{\pi } \over {k+1}}$; $k \in {\bf N}$,$k\neq 0$,
(the so called roots of $1$). Here 
$m=m_{0}+l$, $l \in {\bf Z}$
and we choose $m_{0}$ such that should be contained in one of the $J_{\delta}$
allowed intervals. This does not change the generality.
$m_{0}={{\pi } \over {2s}}-{{\delta } \over {2}}+\epsilon$ 
with $\epsilon $ arbitrary within
$\epsilon \in \biggl (0, (k+1){{2\alpha} \over {\pi}}-k \biggr )\in (0,1)$
such that $m\in \bigcup_{k} (J_{\delta }\cup J_{\Delta })$ only.
Here $k$ fixes the number of states in each of 
the intervals $J_{\Delta }$ and $\Delta\in (k-1,k)$. 
In the interval $\delta $ we have always one single
state. We note that these series \underline{ (2a)},
 labeled by $c$ and $\epsilon $,
 are similar with the complementary
series of $su(1,1)$, because the eigenvalue of the Casimir takes
(continuously) any value
in the interval 
$(c_{1},c_{0})=\biggl ({1 \over {4\sin ^{2}(s/2)}},
{1 \over {\sin ^{2}s}}\biggr )$.
A peculiarity of these series is that the values of the numbers
$N_{\pm }^{cm}$ have exact only $N+1$ distinct values, which repeat after
each sequence of states of the form $|c,m>,...|c,m+k+1>$, or, in terms of
intervals, repeats after each $J_{\delta }\cup J_{\Delta }$ sequence.
Such an example is plotted with full dots in Fig.1a, for $k=1$,
in the range $c_{1}<c<c_{0}$.

Another set of infinite-dimesional unitary representations
 can be constructed for
$s={n \over {k\pi }}$, $n,k\in {\bf Z}, n<k$, for the same range of $c$.
The procedure is identical with the previous one with the difference that
the periodical structure of the $m$'s will cover more than one $J_{\delta
}\cup J_{\Delta }$ sequence. This fact is easy to prove if we project 
the periodical sequence $J_{\delta }, J_{f}, J_{\Delta}, J_{f}$ on a circle
of radius $1/2s$ and find that always exists an integer number $N$ such that 
for
an appropriate $m_{0}$ the label $m_{0}+j$, $j=0,1,...N$ generates a
 periodical sequence 
on this circle. Since $N$ is finite we can always find a set of $m_{0},s$ and
$c$ such that the forbiden intervals do not contain any state of this
sequence. When $q$ is, in general, a root of $1$ such representations are
called "cyclic" [29,30,32,33].

We do not analyse in detail the situation when $s=\xi \pi $ with $\xi $  
irrational.
Since $\biggl [ m+{1 \over 2} \biggr ]$ is randomly distributed
in [0,1] we have for $c>c_{0}$ the continuous representations of class
\underline{(1)}.
The period of $m$ (i.e.1)and of $\biggl [ m\pm {1 \over 2} \biggr ] $
are incommensurate.
We only mention that in the case of the upper limit, $c={1 \over {\sin ^{2} s}}$
we obtain  representations with a special behaviour in the limit
$m \rightarrow \infty$, ("strange series"). 
In order to have
infinite-dimensional unitary representations,
 we need, for any $m$, the condition $N_{\pm
}^{cm}\neq 0$, i.e. $m\neq m_{f}={{\pi} \over {s}}\biggl (k+{1 \over 2}
\biggr )\pm {1 \over 2}$. 
We can write the allowed values of $m$ like $m_{l}=m_{0}+l$, $l\in Z$,
$m_{0}$ arbitrary.
Due to the irrationality of $s/\pi $ we can never have $m_{l}=m_{f}$, for any
finite $l$ and $k$. However, in the limit $l,k \rightarrow \infty $ the
sequences of $m_{l}$ and $m_{f}$ could have the same "accumulation point" 
in the sens that for some sub-sequences of $m_{l}$ and $m_{f}$ we have
$\lim_{l,k \to \infty }{{m_{l}} \over {m_{f}}}=1$. This gives the "strange"
behaviour of such series at infinity and it is an aspect which could be
correlated with the problem of random potential theory.

\vfill

{\it b. Finite dimensional unitary representations of dimension N+1}.

(c bounded and discrete, m bounded)

These unitary representations exist for any $s\in (0,\pi )$. The values of $m$ are
only in the ${\bar J}_{\Delta } $ domains and $m=m_{0}, m_{0}+1,...,m_{0}+N$.
The value $m_{0}$, which fixes the begining of each series, is given by:
$m_{0}={{\pi } \over {s}}-{{\Delta } \over {2}}+{{k\pi} \over {s}}$
with $k\in {\bf Z}$, i.e. different non-equivalent representations for
different $k$'s. The Casimir eigenvalues are
 $c=\biggl [{{N+1} \over {2}}\biggr ]^{2}$
with $N$ given by the condition $\cos ^{2}{s \over 2}<\sin ^{2}s\biggl (
{{N+1} \over 2}\biggr )$. An example for $N=4$ is plotted in Fig.1a, with
stars.

{\it c. Finite dimensional unitary representations of dimension 1}.

These representations belong to the situations when $m$ is 
only in the $J_{\delta }$
intervals and $c=c_{2}=
{1 \over {4\sin ^{2}(s/2)}}$. In this case $m=m_{0}$ only, with $m_{0}=
{{\pi } \over {s}}+{{k\pi} \over {s}}$, $k\in {\bf Z}, \delta =0$
and $s$ takes any value in $(0,\pi )$.

\vskip 0.7truecm
\underline{3. Discrete series,
only finite-dimensional unitary representations, Fig.1b}

(c discrete, bounded and m bounded)
\vskip 0.7truecm
This is the case when $c_{2}=
{1 \over {4\cos ^{2}(s/2)}}<c<c_{1}={1 \over {4\sin ^{2}
(s/2)}}$, and only finite dimensional series exist, due to the fact that
the forbiden intervals for $m$ are larger than 1. This case is valid for
$s\in (0,\pi /2)$. At $s=\pi /2$ the allowed domain vanishes
and for $s\in (\pi /2, \pi )$
we obtain the same numbers as for $s\in (0,\pi /2)$ but the ends of the
allowed interval for $m$ interchange their expressions due to the symmetry of
the equations.
The dimension of each unitary representation in this class is given 
by $N\in \biggl (
0, {{\pi } \over s}-2\biggr )$, depending on the $c$ value in this range.
We have $m\in {\bar J}_{\Delta }=\biggl [{{\pi } \over {s}}+k{{\pi } \over {s}}-
{{\Delta } \over {2}},{{\pi } \over {s}}+k{{\pi } \over {s}}+
{{\Delta } \over {2}}
 \biggr ]$. In each of these cases we can extra-label the number $m$ by $m^{(l)}=
{{\pi } \over {s}}+k{{\pi } \over {s}}-
{{\Delta } \over {2}}+l $ for $l=0,1,...,N$ ($k\in {\bf Z}$ labels the
representation). We have also $c=\biggl [{{N+1} \over 2} \biggr ]^{2}$
or equivalently $c= \biggl [j+{1 \over 2} \biggr ]^{2}$ if we denote
$N/2=j$ and consequently the dimension of the unitary representations are $N+1=2j+1$
with $N_{max}={{\pi} \over {s_{min}}}-2$, similar with the results 
presented in [33]. An example with $N=5$ is ploted in Fig.1b.
\vskip 0.7truecm
\underline{Comments and examples}
\vskip 0.7truecm
The above presented clasification is in agreement with the traditional
clasification in the existing literature, i.e.: spliting of the
unitary representations in the heighest weight (of type 
\underline{(2b)},
\underline{(2c)} and
\underline{(3)} in our case) and cyclic ones (of type 
\underline{(1)} and
\underline{(2a)}), [29,30].

We should like to do some comments and to workout some simple examples
concerning some special values for $s$ and $c$ in these representations.
First we mention that for the last interval of the eigenvalue
$c<{1 \over {\cos ^{2}(s/2)}}$ there are no unirreps.
Second, we mention that the mixed series \underline{(2)} and 
the discrete series \underline{(3)} coincide with those given in [30,33]
 for the case
when $q$ is a root of the unity.

For both the mixed \underline{(2)} and discrete \underline{(3)} series 
presented above, we have the action
\begin{eqnarray}
J_{\pm}|c,m> & =& \sqrt{c-\biggl [ m\pm {1 \over 2}\biggr ]^{2}
}|c,m\pm 1>, \hfill
\end{eqnarray}
according to
the corresponding unirreps in [5], up to an undetermined phase factor.

We note that though we have deformed a compact
Lie algebra $su(2)$, the above q-deformation allows the existence of
continuous series for $su_{q}(2)$, equivalent with $C_{x}^{0}$
principal series of $su(1,1)$, up to a redefinition of the Casimir,
$C \rightarrow C+{1 \over 4}-{1 \over {\sin ^{2}s}}$.
Different from the notations of the above cited papers, we label the
representations with the $c$ number instead of $j$. We obtained continuous
series,
not only for the $s\rightarrow \pi$ ($q=-1$) case like in [34] 
but also for other
values of the deformation parameter $s$ (class \underline{(1)} and 
\underline{(2a)}). 

We want to note also the existence
of some special deformations. In the limit $s\rightarrow \pi $ we have
\begin{eqnarray}
[2J_{z}]|cm> & =& 
\cases{-2m|cm>  & $ m\in {\bf Z} $
\cr
2m|cm> & $m\in {\bf Z}+{1 \over 2}$
\cr
\infty & $ \mbox{in the  rest} $
}.
\hfill
\end{eqnarray}
In this last case we can explain what we mean by the similarity between
 $su_{e^{i\pi }}(2)$ and $ su(1,1)$,
for $m\in Z$,
through eq.(7). 
For $s \rightarrow \pi $ and $m\in Z$ the unitary representations
\underline{(2a)} of $su_{e^{is}}(2)$ are equivalent with the $C_{x}^{0}$
unirreps of $su(1,1)$: in both these representations $m$ has the same
range and the commutators are identical on the basis $|cZ>$.
The single apparent problem consist in the fact that the Casimir
operator has poles. In order to obtain the equivalence we can introduce a
manifold of smooth functions $\sigma :[0,\pi ] \rightarrow R$ such that $\sigma
(\pi )=0$, $\biggl | {{d\sigma } \over {ds}} \biggr |(\pi )<\infty $.
This manifold
 is continuous with respect to a certain parameter which can label
these functions, e.g. $\alpha ={{d\sigma } \over {ds}}(s=\pi ) $.
In this way, for any $s$ we can choose a certain element of this class, $
\sigma _{\alpha }(s)$, such that $|j(\sigma )m>=|\alpha m>$ and the functions
$\sigma _{\alpha }(s)$ are choosen such that the limit $s\rightarrow \pi $
exists for $C$ on this basis.

For $s \rightarrow \pi / (k+1)$ , $k\in {\bf N}$, $k\neq 0$, the action of $[2J_{z}]$ is
$$[2J_{z}]|c,m>={{\sin 
{{(2{\tilde m}+1)\pi } \over {(k+1)}}
} \over {\sin
{{\pi } \over {k+1}}}}|c,m\pm 1>,
$$
where ${\tilde m}\equiv m_{(k)}$ is $m$ modulo $k$. Consequently there are
only $k$ distinct eigenvalues for $[2J_{z}]$, and also for $J_{\pm }$.
The states $|c,m>_{s=\pi/(k+1)}$ have infinite degeneration for each value
of $m$. We present graphicaly such an example in Fig.2, where the first 9
eigenvalues
 of $[2J_{z}]$ are ploted against $s$.
In the case when $s=\pi /2$  we have
\[ [2J_{z}]|cm> = \left\{ \begin{array} {ll}
\pm 1|cm>  & \mbox{ $m\in {\bf Z}+{1 \over 2} $ \quad (quasi-spin-like
 representations)} \\
0 & \mbox{$ m\in {\bf Z}$  \quad \quad \quad
(e(2)-like commutators) }   
\end{array}
\right. \]
We note that in the second case of the above relation,
the commutators of the resulting $su_{i}(2)$ algebra, i.e.
eq.(1) and
\begin{eqnarray}
[J_{+},J_{-}]=0,
\end{eqnarray}
provide an algebra homeomorphism with the Euclidean Lie
algebra $e(2)$, [35], on the subspace of the unitary representations 
$|c, {\bf Z}>$ only.
We mention that previously  we have used the word "similarity" 
instead of
algebra homeomorphism, because in the q-deformed cases the action of
$[2J_{z}]$ is identical with the action of $2J_{z}$, only on  same 
certain subspaces of the representations,
for
the corresponding undeformed Lie algebras (i.e. $su(1,1)$ and $e(2)$ in the
above examples). Evidently a homeomorphism involving $2J_{z}$ and its
q-deformation could be established at the level of the universal covering
algebra, but such an analysis is beyound
the aim of this paper. Other examples of
such homeomorphisms between $su_{e^{s}}(2)$ and other non-linear
algebras are investigated in [25,27].
More, as a support
of the above comments, we can reobtain, in our case, the Maekawa's 
equations
for the Casimir, [9], for $q=e^{2is}$
\begin{eqnarray}
C & = & \cos s[J_{z}]^{2}+{1 \over 2}(J_{+}J_{-}+J_{-}J_{+}) \hfill
\end{eqnarray}
$$=\cos s [J_{z}]^{2}+J_{x}^{2}+J_{y}^{2}+{1 \over {4\cos ^{2}(s/2)}}.$$
We note that in the limit $s \rightarrow \pi /2$ we obtain
\begin{eqnarray}
C_{s=\pi /2} & \rightarrow & J_{x}^{2}+J_{y}^{2}+{1 \over 2}, \hfill 
\end{eqnarray}
which gives the Casimir of the $e(2)$ Euclidean Lie algebra, in
agreement with the consequences of eq.(18) (see also [36] for $e_{q}(2)$).
In the case $s \rightarrow \pi $ we obtain from eq.(19)
$$C=-J_{z}^{2}+{1 \over 2}(J_{+}J_{-}+J_{-}J_{+}),$$
but again only on the sub-space of states $|c, {\bf Z}>$, 
which, together
with the first line in the RHS of eq.(17) and with the comments in the
construction of the
unirreps of class \underline{(1)} and \underline{(2)}, 
support the asumption for the similarity
of $su_{e^{i\pi }}(2)$ with $su(1,1)$.
This is in agreement with the conclusions
in the references [22,34], too.
We would like to note that such a behaviour could be a good start for an
analysis of the transition from bounded to unbounded states
of a system having q-dynamical symmetry, in the AST frame.
These features can provide also a new possibility for the introduction
of the Euclidean connection, as a natural extension of the initial
undeformed symmetries, [19,20,28,34,36].

However we would like to mention another consequence of eq.(19),
 which supports the
comparison between $su_{e^{is}}(2)$ and $su(1,1)$.
According to the geometrical interpretation of Alhassid {\it et al}
in the second reference in [28], the surfaces of constant eigenvalues of the
Casimir operator, described in a formal 3-dimensional vectorial
space of $(J_{x},J_{y},J_{z})$ , have (for $su_{e^{is}}(2)$) the form
\begin{eqnarray}
C& = & \cos s {{\sin ^{2}(sJ_{z})} \over {\sin ^{2} s}}+J_{x}^{2}
+J_{y}^{2}=const.
\end{eqnarray}
Loosely speaking, these formal "surfaces" are a sort of deformations 
of the sphere
associated with the geometry of the Lie group $SU(2)$. For
$c$ and $s$ in the range of the infinite dimensional unirreps
\underline{(2a)}, they form an infinite countable
disconnected set of identical closed sufaces. In the infinite dimensional
unirreps case \underline{(1)}, with the increasing of $s$,
these surfaces become connected together into one single connected surface,
i.e. an open surface, homotopic with the hyperboloid surface associated with
the non-compact Lie group $SU(1,1)$, [28]. Of course, in our approach
 one can not give a 
rigurous group theoretical or algebraic 
interpretation to these algebraic surfaces because, in the q-deformed case,
there is no associated Lie group manifold, or at least in this approach.
We illustrate a two-dimensional section (i.e. in the $(J_{z},J_{x})$-plane) 
of such "surfaces", for the compact and non-compact limiting cases
in Fig.3. In other words the topology of these "surfaces" is controled
by the parameter $s$ (i.e. $q$) when ranges from $0$ to $\pi $.

\vskip 1.5truecm
\section{Generalised real deformation for $su(2)$ as a Hopf 
algebra }
\vskip 1truecm
In this Section we intend to find other connections of $su(2)$
with different nonlinear algebras and to
generalise the exponential deformation of $su_{q}(2)$
to a polynomial one and to find a
possible differential realisation of it.
The next possibility which we sugest for deforming the bounded 
unitary representations 
of $su(2)$ into 
infinite-dimensional ones is to keep the deformation parameter $q\in {\bf
R}$ and to generalise the deformation function, i.e. the RHS of eq.(2).
Several succesfull attempts were made in these direction [25-27,37,38]
and new kinds of nonlinear algebras, together with their unirreps
were constructed (most of them having exponential-like spectra for
$[2J_{z}]$).
In the following we apply the deformation introduced by Ludu et al. in [26],
mainly for the fact that it can be endowed with a real-hopf algebra structure.

Starting from $su(2)$ we construct a deformed (nonlinear) algebra $L_{g,q}
\in U(su(2))$ (the universal covering algebra), q real, 
generated by the operators $J_{z}=J_{z}^{\dagger}$, $J_{+}$, $J_{-}=
J_{+}^{\dagger}$,
the identity I, $g$ and $g^{-1}$. The last two elements
can also be considered as  two
holomorphic real functions $g(J_{z},q),g^{-1}(J_{z},q)$ which belong to the
polynomial ring ${\it C}[J_{z},J^{-1}_{z}]$ and are defined by the relations
\begin{eqnarray}
gg^{-1} & = & g^{-1}g=1,\hfill
\end{eqnarray}
\begin{eqnarray}
g^{2}J_{\pm }g^{-2}&=q^{\pm }J_{\pm},\hfill
\end{eqnarray}
\begin{eqnarray}
[J_{+},J_{-}]&=&2{{g^{2}-g^{-2}} \over {h(q)}},\hfill
\end{eqnarray}
where $h(q)$ is an analytical function of q. It is proved [26] that $L_{g,q}$
satisfies the postulates of the Hopf algebra [29-32] if
\begin{eqnarray}
\Delta (J_{\pm })&=&J_{\pm }\otimes g^{-1}
+g\otimes J_{\pm },\hfill
\end{eqnarray}
\begin{eqnarray}
\Delta(g^{\pm 1})&=&g^{\pm 1}\otimes g^{\pm 1},\hfill
\end{eqnarray}
\begin{eqnarray}
\epsilon (J_{\pm })&=&0, \quad
\hbox{~~~} \epsilon (g^{\pm 1})=1,\hfill
\end{eqnarray}
\begin{eqnarray}
S(J_{\pm })&=&-q^{\mp 1}J_{\pm }, \hfill
\end{eqnarray}
\begin{eqnarray}
S(g^{\pm 1})&=&g^{\mp 1}.\hfill
\end{eqnarray}
The algebra $L_{f,q}$reduces to $su(2)$ in the limit $ q \rightarrow 1$
if the following conditions are fulfilled
\begin{eqnarray}
f(J_{z},1)&=&1,\hfill
\end{eqnarray}
\begin{eqnarray}
{{\partial f} \over {\partial q}}
(J_{z},1)&=&J_{z},\hfill 
\end{eqnarray}
\begin{eqnarray}
h(q)&=&q-q^{-1},\hfill
\end{eqnarray}
where we denote $g^{\pm 2}=f^{\pm 1}$. If $f=q^{J_{z}}$ then $L_{g,q}$ is
homeomorphic with $sl_{q}(2,C)$ for complex q and with $su_{q}(2)$
for real q [1,30]. In this latter case the q-deformed algebra has a
real-Hopf algebra structure defined by the relations (22-32).
\par
The Casimir operator of $L_{g,q}$ is given by [26,27,30]
\begin{eqnarray}
C&=&2^{1/3} \biggl [
{
{
2{(q^{-1/4}g-{q^{1/4}}{g^{-1}})}^{2}
} \over {
{(q-q^{-1})}^{4/3}
{(q^{1/2}-q^{-1/2})}^{2/3}
}}+{{J_{+}J_{-}} \over {
{(q^{1/2}+q^{-1/2})}^{1/3}
}}\biggr ].
\hfill 
\end{eqnarray}
In order to obtain infinite-dimensional representations
for ${\cal L}_{g,q}$ we introduce the definition
\begin{eqnarray}
f(J_{z},q)=1+(q-1)J_{z}+(q-1)^{2}(\alpha J_{z}+b(J_{z})), \hfill
\end{eqnarray}
where $b(J_{z})$ is a positive, bounded analitical function, arbitrary for the
moment. It is easy to verify that the function (34) fulfills the conditions
(30-31). For $q=q_{1}={{\alpha -1} \over {\alpha}}$
we have $f(J_{z},q_{1})=1+(q_{1}-1)^{2}b(J_{z})$.
In order to obtain the unitary
representations for this q-deformed algebra we use the same
condition of unitarity like in Sect.1, i.e. that the matrix elements of the operators
$J_{\pm }J_{\mp }$ should be non-negative. Consequently, on the fundamental
representation $|c,m>$, we obtain the corresponding conditions
\begin{eqnarray}
<J_{+}J_{-}>&=&{1 \over {C_{2}}}\biggl (c-C_{1}({\tilde g}-
{\tilde g}^{-1})^{2}\biggr )
\geq 0, \hfill
\end{eqnarray}
\begin{eqnarray}
<J_{-}J{+}> & = & {1 \over {C_{2}}}\biggl (c-C_{1}({\tilde g}-{\tilde
g}^{-1})^{2}-{{2C_{2}} \over {h(q_{1})}}({\tilde {g}}^{2}-{\tilde {g}}^{-2})
\biggr ) \geq 0,  \hfill
\end{eqnarray}
where
$$C_{1}={{2^{4/3}} \over {(q_{1}-q_{1}^{-1})^{4/3}
(q_{1}^{1/2}-q_{1}^{-1/2})^{2/3}}}\geq 0,$$
$$C_{2}={{2^{1/3}} \over {(q_{1}^{1/2}-q_{1}^{-1/2})^{1/3}}}\geq 0, $$
and
$${\tilde g}=q_{1}^{-1/4}\sqrt{f(J_{z},q_{1})}\geq 0.$$
Both conditions (35,36), reduce to second order inequalities
concerning the function $f$ and ask for the image of this
 function to be bounded in some
regions, depending on $c$ and $q_{1}$. 
For example eq.(35) is equivalent with the condition
\begin{eqnarray}
f & \in & [L_{1}(c,q_{1}), L_{2}(c,q_{1})], \hfill
\end{eqnarray}
with
$$L_{1,2}(c,q_{1})=\sqrt{q_{1}}\biggl ( 
\beta (c,q_{1})\mp \sqrt{ \beta (c,q_{1})
^{2}}-1\biggr ), $$
with $\beta (c,q_{1})=(1+c/2C_{1}(q_{1})$.
These limiting values coincide for $c=0$,
$L_{1}(0,q_{1})=L_{2}(0,q_{1})=\sqrt{q_{1}}$ and tend to
the limits $L_{1} \rightarrow 0$, $L_{2} \rightarrow \infty $, when 
$c \rightarrow \infty $, respectively. 
Consequently the condition (35) is fulfiled for $c$
larger than a certain value, such that the image of $f$ should be included
in $(L_{1},L_{2})$. The second condition, eq.(36), can also be 
written in the form
\begin{eqnarray}
f\in [l_{1}(c,q_{1}),l_{2}(c,q_{1})], \hfill
\end{eqnarray}
with
$$l_{1,2}=\sqrt{q_{1}}(\beta ^{2}\alpha _{+}\mp\sqrt{\beta ^{2}\alpha
_{+}^{2}-{{\alpha _{+}} \over {\alpha _{-}}}
},$$
with $\alpha _{\pm }^{-1}=1\pm 2C_{2}/C_{1}(q_{1})h(q_{1})$.
The analysis of eqs.(37,38) shows that for $c$ larger
than a certain positive value, $c_{min}(q_{1})$ ,
there are always nontrivial solutions for $f$, i.e.
$l_{2}>L_{2}>l_{1}>L_{1}$
for $c>c_{min}(q_{1})$.
We denote the boundaries of the allowed domain of $f$ 
with $f\in (f_{min}(c,q_{1}),f_{max}(c,q_{1}))$.
This interval exists for $c\geq c_{min}(q_{1})>0$ and its lenght 
increases with the
increasing of $c$. This finaly proves that one can chose
the functions $f$ in eq.(34) such that their images are 
included in $(f_{min},f_{max})$.
Consequently, the unitarity conditions (35,36) are fulfiled for all the values
of $f(J_{z},q)$ in eq.(34), i.e. for all $m$ eigenvalues
and $q=q_{1}$. 
Since there are no limitations concerning
$m \in Z$ or ${1 \over 2}Z$ this proves the existence
of infinite-dimensional unitary representations 
of ${\cal L}_{\tilde g}, q_{1}$ described by $m\in {\bf Z}$ or ${1 \over
2}{\bf Z}$ (unbounded from below and above), $c\geq c_{min}(q_{1})$, 
(continuous, bounded from
below and unbounded from above) and $g^{2}$ given by eq.(34).
In the following, we give a simple example of such a q-deformation.
We can choose for the function $f$ in eq.(34) the form 
\begin{eqnarray}
f(J_{z},q_{1})& = & {{f_{max}(q-{1})-f_{min}(q_{1})} \over 
{\cosh ^{2}(J_{z})}}+f_{min}(q_{1}). \hfill
\end{eqnarray}
It is now easy to calculate the corresponding spectrum of $[2J_{z}]$,
 for this deformation.
This spectrum is discrete, infinite dimensional,
contains one accumulation point and it is bounded from
above and below. 
Such spectra are very usefull in physical applications
(atomic or molecular physics) because they can generalise
atomic or molecular spectra where the distance between two consecutive lines
continuously decreases. The Hidrogen atom has, as an example, such
a kind of spectrum with its accumulation point at its superior
bound. There is also interest for such spectra in conformal field
theories.

We would like to do some comments at the end of Sections 2 and 3.
We have shown that starting from the $su(2)$ Lie algebra symmetry, we
could obtain deformed algebras (through q-deformation with $q=e^{is}$
in Sect.2 or through a general real q-deformation of the commutators 
in Sect.3), having infinite-dimensional unitary representations, some of
them being irreducible too.
The transition between the bounded unirreps ${\cal D}^{j}$ of $su(2)$
to the unitary representations of the corresponding
deformation is made by the continuous variation of the parameter of
deformation, in both cases.
Consequently, physical systems having dynamical symmetry close to that of
$su(2)$ could be described in this formalism. We expect especially to obtain
good results in the investigations of  scattering processes, where
the Hamiltonian could be expressed in terms of the deformed operators,
like $[2J_{z}]$ or $C$. In this case, the parameter of deformation could play
the role of a coupling parameter of the system, which allows the smooth
transition from bounded to free states. A more complete picture could be
obtained if one can construct differential representations of such deformed
structures. In this way a Schr\"odinger picture can be obtained together
with informations about the corresponding potentials and, the
interactions. This is the task, for the above discussed case, in the
following Section.

\vskip 1.5truecm
\section{One-dimensional potential realisations
 for the complex deformation $su_{e^{is}}(2)$ }
\vskip 1truecm

In this Section we introduce 
a differential operator representation of $su_{e^{is}}(2)$, $s\in R$,
acting on a two-dimensional space, in order to obtain a 1-dimensional 
Schr\"odinger equation and, consequently, equivalent potential picture for
different values of $s$ and $m$.
We strictly follow the procedure given in [28], but
we use a sort of "q-deformation" of the derivative with respect
to one of the coordinates, as was early suggested by Biedenharn [4] and
Macfarlane [5]. 
In this way one could use the $su_{e^{is}}(2)$
algebra as a dynamical symmetry close to some
physical applications. We are interested especially 
in the classes of potentials 
which are connected with the
eigenfunctions of the
infinite-dimensional unitary representations obtained in Section 2.
\par
We introduce a 2-dimensional coordinate space, parametrised with
the polar coordinates $\phi
\in [0,\pi )$, $r \in R$. The differential realisation
is given by
\begin{eqnarray}
J_{\pm}&=&e^{\pm i\phi } \Biggl ( \pm
\partial_{r}+f_{1}(r)
 \Biggl ( {1 \over 2}
{{[2i\partial_{\phi}]}}
 \Biggr ) +f_{2}(r) \Biggr ), \hfill
\end{eqnarray}
\begin{eqnarray}
J_{z}&=&-i\partial_{\phi}, \hfill
\end{eqnarray}
where the $ \pm$ signs appear in order to satisfy the hermiticity
conditions for $J_{\pm },J_{z}$ and
 $f_{1},f_{2}$ are real smooth functions of $r $ 
which have to be determined 
 in order to re-obtain the commutators of $su_{q}(2)$. The operator 
 $[2i{\partial }_{\phi }]=-[2J_{z}]$
is a sort of q-deformation
of the derivative  with respect to the $\phi$ coordinate. 
This q-deformation is defined by the formal 
action of the exponential 
operator (Taylor series) on
exponential functions. By introducing eq.(41) in 
eq.(9)
we have the action
\begin{eqnarray}
[2i{\partial }_{\phi }]e^{im\phi }&=&-[2m]e^{im\phi }, \hfill
\end{eqnarray}
for $q=e^{is}$, $s$ real, $[m] \in {\bf R}$. Consequently
$[2i{\partial }_{\phi }]$ commutes with $J_{z}$ and $[J_{z}]^{\dag }=
[J_{z}]$.
\par
In the following, we deduce some properties of the operator 
$[2i{\partial }_{\phi }]$. We note that this operator 
is self-adjoint on the set of 
the functions of the
form
$|jm>=\Psi_{jm}(r,\phi )=R_{jm}(r)e^{im\phi }$ 
 We introduce the auxiliar operators
\begin{eqnarray}
{\Omega }_{\pm }&=&e^{\pm i\phi }[2i{\partial }_{\phi }]
e^{\mp i\phi }, \hfill
\end{eqnarray}
which act through eq.(42) like
\begin{eqnarray}
{\Omega }_{\pm }e^{im\phi }&=&-[2(m\mp 1)]e^{im\phi }.
\hfill
\end{eqnarray}
We need some auxiliar relations, defined for any two generic 
quantities $x_{1},x_{2}$, in the form
\begin{eqnarray}
[x_{1}-x_{2}]+(-1)^{k}[x_{1}+x_{2}]&=&
(-1)^{k}{\eta }^{2}[x_{3-k}] {{\biggl [{{x_{k}} \over 2}\biggr ] }^{2}}
+2(-1)^{k}[x_{3-k}], \hfill
\end{eqnarray}
and
\begin{eqnarray}
[x_{1}+x_{2}]+(-1)^{k}[x_{1}]&=&{\eta ^{2}} \biggl [
{{x_{2}} \over {2}}\biggr ]^{k} {\biggl [{1 \over {3-k}} 
\biggl (x_{1}+
{{x_{2}} \over 2}\biggr )
\biggr ] }^{3-k}
+ 2\biggl [ {{x_{2}} \over 2}+(k-1)x_{1} \biggr ], \hfill
\end{eqnarray}
with $k=1,2$ and $\eta =2i\cdot sin(s)$.
In order to realise the $su_{e^{is}}(2)$ algebra withthe operators
$J_{\pm }, J_{z}$, given in eqs.(40-41), these operators must satisfy the
commutators of this deformed algebra, i.e. eqs.(1,2).
By using eqs.(45-48) it is easy to check that 
these operators satisfy eq.(1).
For eq.(2) we have $(f_{i}^{'}=df_{i}/d\rho )$
\begin{eqnarray}
[J_{+},J_{-}]e^{im\phi }&=&
-[2m]\biggl (f_{1}^{'}+{{f_{1}^{2}} \over 2}
\biggr )e^{im\phi }+{\bf O}e^{im\phi }, \hfill
\end{eqnarray}
where the operator ${\bf O}$ acts on $e^{im\phi }$ in the form
\begin{eqnarray}
{\bf O}&=&\biggl [ [2]f_{1}f_{2}
\biggl ( 1+{{{\eta }^{2}{[m]}^{2}} \over 2}\biggr ) -{{[2]} \over 4}f_{1}^{2}{\eta }^{2}
[2m][m]^{2}+2f_{2}^{'}+{{f_{1}} \over 2}{\eta }^{2}[2m]{\partial }_{\rho }
 \biggr ]. \hfill
\end{eqnarray}
Consequently, by identifying eq.(2) with eq.(47),
we obtain two conditions
\begin{eqnarray}
f_{1}^{'}&=& -1-f_{1}^{2}cos(s), \hfill
\end{eqnarray}
\begin{eqnarray}
{\bf O}R_{cm}e^{im\phi }& = &0. \hfill
\end{eqnarray}
The action of the Casimir operator of 
$su_{e^{is}}(2)$, eq.(3), on the basis $R_{cm}(r)e^{im\phi }$, reduces to
$$\biggl ( -{\partial }_{r}^{2}+ {{f_{1}} \over 2}
\biggl (2+{\eta }^{2}\biggl [m-{1 \over 2}\biggr ]^{2}\biggr )
{\partial }_{r} +[2m][2(m-1)]{{f_{1}^{2}} \over 4}- $$

$$
{{f_{1}f_{2}}
\over 2}
(2+{\eta }^{2})[2m-1]-{{f_{1}^{'}} \over 2}[2m]+f_{2}^{2}+f_{2}^{'}+
[m]^{2}+$$
\begin{eqnarray}
\biggl [ m-{1 \over 2} \biggr ]
\biggr )R_{cm}(r)&=&cR_{cm}(r), \hfill
\end{eqnarray}
for any $m$, in the corresponding representation labeled by $c$.

The first condition, eq.(49), is similar
with the corresponding condition usualy used in differential realisations
of $su(2)$ or $su(1,1)$ in AST [28,40], excepting the deformed coefficients
like $[2]/2$, instead of $1$. This equation has unique solutions,
independent of the quantum numbers.
There exist
three different cases, accordingly to the range of $s$.
For $cos(s)<0$ (which includes the limiting case $s \rightarrow \pi $
of $su(1,1)$) we have
\begin{eqnarray}
f_{1}(r)&=&\pm {1 \over {\sqrt{-cos(s)}}} \cdot
\cases{tanh(\mp \sqrt{-cos(s)}r +d) & $ $
\cr
1 & $ $
},
\hfill
\end{eqnarray}
where $d$ is a constant of integration.
For $cos(s)>0$ (which contains the limiting case $s \rightarrow 0$ of $su(2)$)
we have only one  solution
\begin{eqnarray}
f_{1}(r)&=&\pm {1 \over {\sqrt{cos(s)}}} 
tg(\mp \sqrt{cos(s)}r+d), 
\hfill
\end{eqnarray}
and no constant real solution.
Finally, for $cos(s)=0$ we have
\begin{eqnarray}
f_{1}&=&-r+d. \hfill
\end{eqnarray}
In order to have the continuity of $f_{1}$ between
 the above solutions, with respect to
$s\in [0,\pi ]$, we need to fix $d=0$.
We note that the solutions in eqs.(52-54) are similar with those 
obtained in the
undeformed case, excepting the scalling term $\sqrt {\pm cos(s)}$.

The second condition, eq.(50), is more special due to the fact that 
the operator
${\bf O}$ from eq.(48) contains the derivative with respect 
to $r$. This a special
effect due to the deformation (does not occure in the undeformed case),
i.e. due to the terms containing $\eta \neq 0$.
This special dependence
introduce a coupling between the function $f_{2}$ and the eigenstates
$R_{cm}$.
Consequently,
we can see that, due to the deformation, a different situation compared with
the undeformed case appears:
one has the possibility of obtaining other
different
potentials from the $su_{e^{is}}(2)$ dynamical symmetry, 
in comparison with the undeformed
case. On the other hand it is more difficult to realise the eigenproblem for
the Casimir operator into a Schr\"odinger equation, due to the fact that
there are two coupled differential equations for the eigenstates, 
i.e. eqs.(50,51).
\vskip 1truecm
\underline{Limiting cases}
\vskip 0.5truecm
We analyse this differential realisation for three (simpler) limiting cases,
i.e. when eqs.(50) and (51) become decoupled due to the canceling of
the last term in the RHS of eq.(48). This happens for
$s\simeq \pi$ (any $m$), $s\simeq \pi /2$ ($m\in Z$) and 
$s\simeq 0$ (any $m$) with 
$s\in (0,\pi)$.
We supose that for all
these situations the coefficients $\eta ^{2}$ and/or $[2m]$ are
 enough small to be
neglected.
Consequently, in these approximations, eq.(50) can be integrated independent
of $R_{cm}$ and we obtain an approximate differential equation for $f_{2}$
\begin{eqnarray}
{{[2]} \over 2}f_{1}f_{2}&=&-f_{2}^{'}. \hfill
\end{eqnarray}
By integrating this equation with $f_{1}$ given by
eqs.(52) we have, for the first limiting case, ($s\simeq 0$)
\begin{eqnarray}
f_{2}(r)&=&
\cases{
{{F_{1}} \over {cosh(\mp \sqrt{-cos(s)}r +d1)}} & $ $
\cr
F_{2}e^{\pm \sqrt{-cos(s)}r} & $ $
}.
\hfill
\end{eqnarray}
For $s\simeq 0$ limiting case ($cos(s)>0$) we have
\begin{eqnarray}
f_{2}(r)& = & F_{3}(cos (\mp \sqrt{cos(s)}r+d2)). \hfill
\end{eqnarray}
Finally, in the case $s\simeq \pi /2$ we have
$$f_{2}(r)=F_{4}=const.$$
The constants of integration $F_{1-4}$ and $d1, d2$
can be choosen such that $f_{2}$ should be continuous with respect to $s$.
With the above solutions for the functions $f_{1}$ and $f_{2}$ we can
construct the differential form of the Casimir operator, close to the
limiting cases $s\simeq 0$ and $s\simeq \pi$
(i.e. in the approximation $\eta ^{2} << 1$)
$$
Ce^{im\phi }=\biggl (
-{\partial }^{2}_{r}-f_{1}{\partial}_{r}
+[2m][2m-2] {{f_{1}^{2}} \over 4}-f_{1}f_{2}[2m-1]
$$
\begin{eqnarray}
-{{f_{1}^{'}} \over 2}[2m]+f_{2}^{2}+f_{2}^{'}+[m]^{2}
+\biggl [ m- {1 \over 2} \biggr ]^{2}
\biggr ) e^{im\phi },
\end{eqnarray}
and close to the intermediate limiting case $s\simeq \pi /2$, $m\in Z$ 
$$
Ce^{im\phi }= -{\partial }^{2}_{r}-{{f_{1}} \over 2}(2+(-1)^{m+1})
{\partial }_{r}+(-1)^{m}{{f_{1}f_{2}} \over 2}
(2+{\eta }^{2})+$$
\begin{eqnarray}
f_{2}^{2}&+&f_{2}^{'}+[m]^{2}+1)e^{im\phi }. \hfill
\end{eqnarray}
Both eqs.(58,59) can be reduced to a 1-dimensional Schr\"odinger-like
equation
\begin{eqnarray}
(-{\partial }^{2}_{r}+V(r;m,s)){\Psi }_{cm}(r)&=&c{\Psi }_{cm}(r),
\end{eqnarray}
${\Psi }_{cm}(r,\phi )=R_{cm}(r)e^{im\phi }$ 
by using the substitution: $R_{cm}(r)\rightarrow a(r)R_{cm}(r)$.
For example, in the case of eq.(58) we choose
\begin{eqnarray}
a(r)=a_{0}exp\biggl (-\int f_{1}(r)dr \biggr ),
\end{eqnarray}
and we obtain a potential dependent of $r$ and of the parameters $m$ and $s$
$$
V(r;m,s)=-{{a^{''}} \over a}+[2m][2m-2]{{f_{1}^{2}} \over 4}
-f_{1}f_{2}[2m-1]-
$$
\begin{eqnarray}
{{f_{1}^{''}} \over 2}[2m]&+&f_{2}^{2}+f_{2}^{'}+[m]^{2}+
\biggl [ m-{1 \over 2}  \biggr ]^{2},
\end{eqnarray}
with $f_{1},f_{2}$ given by eqs.(52-54) and (56,57), respectively.
Eq.(62), gives a posible equivalent 1-dimensional 
potential picture (Schr\"odinger) for the limiting cases $s\simeq 0$
and $s\simeq \pi$ and close to them.
In the case $s\simeq \pi /2$, $m \in Z$, one can use the same 
type of substitution and
$a(r)$ differs only by a constant factor in the exponential, in
eq.(61).

For exemplifying we present in Figs.4-6 some potential shapes
drawn for few different limiting values of $s$.

In Fig.4a the potential $V(r;m,s)$ given by eq.(62) is presented in the case
$cos(s)>0$, for different values of $s\in (0,pi /2)$, $s\simeq 0$, $s\simeq
\pi /2$ and $m=1$. The potential consists in an infinite 
series of periodic copies of
well potentials valleys of infinite depth. The aspect of $V$ does not
modifies essentially when $s$ modifies from zero to $\pi /2$. In Fig.4b
the same potential is presented for fixed $s=0.25$ and different $m=1,3/2,
2,5/2,3$ and $7/2$. We note that the modification of $m$, below or above
a limiting value (situated next to $m=2$), drastically transforms
the aspect of $V$: from a countable succesion of infinite-depth potential
wells (at small $m$'s) into a similar sequence of positive poles. 
The shapes for $m=1,3/2$ and $2$ are similar with the potential of a
quantum particle (electron) in an atomic crystal.
One knows that in this case the spectrum is spread out in bands of closely
spaced levels for $c$ (in between the wells) and continuous for $c$ above 
their common tops. 
In the present case the bands are generated
by the discrete representations of class  \underline{(1)}
and \underline{(2)}.
The limiting value for the transition bands/continuous increases
with $s$, as predicts the analysis of these representations, Fig.1a.
For $m=2,5/2,3$ and $7$ we obtain a sum of positive potentials which are
periodicaly distributed on the $r$ axis.
This results from the fact that in the expression, eq.(62) of $V$ there
appear only the functions $f_{1,2}(cos(s)r)$, given by eqs.(53,57).
All these functions have the same trigonometric structure, with the same
period.
This is a consequence of the fact that, due to the decoupling of $f_{2}$
and $R_{cm}$, $f_{2}$ can be directly obtained from eq.(55) and it depends
only on $f_{1}$.
Thus, the resulting potential is only a sum of potentials with the same
space-period. Consequently, the spectrum is absolut continuous [41],
the eigenvalues are of Bloch-type states in which the paticle is found with
infinitesimal probability in every finite region of the axis.
The q-deformation, in this case, ($s$ close to $0$ or to $\pi /2$)
does not bring but quantitative modifications against the undeformed case.
Here $s$ can be regarded as a fitt parameter, only.

In Fig.5 we present some potential shapes in the range $cos(s)<0$,
close to the limits $s\simeq \pi /2$ and $s\simeq \pi $, constructed by
introducing the first solutions of $f_{1,2}$ from eqs.(52,56) in eq.(62).
The resulting potentials are no more periodic and they have exponential
(hyperbolic) behaviour.
In Fig.5a we show the potential for a fixed $m=1$ and for different
$s$. In this figure the potential shapes with $s=3.05$ and $3$ 
($\simeq \pi $) are
drawn translated along the vertical axis, in order to have a better pictorial
view. We note that in the both limits, the potentials modify drastically:
closed to $s=3$, i.e. $s$ close to an integer, behave as pure
repulsive for $r>0$ and, close to $s\simeq \pi $ the potential
 always has a pocket,
close to $r\simeq 0$. The first examples are similar with 
P\"oschl-Teller-like potentials and the latter cases give shapes similar with
a negative signed harmonic oscillator potential.

In Fig.5b the same potentials as in Fig.5a are presented, for fixed $s=3$
($\simeq \pi $)
and different $m$'s. For example, for $m=1$ the potential is 
practicaly symmetric
in $r$, P\"oschl-Teller-like. For larger values of $m$ the potential become
rather antisymmetric in $r$ and consists in a bounded valley (close to
 $r=-1$ for
$m=3$ and close to $r=0.5$ for $m=5/2$)
followed by a shoulder.
In the $s\rightarrow \pi $ limit these potentials coincide with a 
P\"oschl-Teller one, like in the undeformed case.

By introducing the second solutions of $f_{1,2}$ from eqs.(52,56) in
eq.(62),
we get now Morse-like  deformed potentials.
In Fig.6a we present such shapes for $s\simeq \pi /2$ and fixed $m=1$.
Here, all the potentials are also drawn translated along the vertical axis
such that they have the same asymptotic value for $r\rightarrow -\infty$,
i.e. $V(s=3.04)\rightarrow V(s=3.04)-94.8$, etc.
For values of the Casimir eigenvalues samller than the asymptotic 
limit, the spectrum is discrete ($c_{2}<c<c_{0}$), for $c>c_{0}$ the
spectrum is continuous and for $c$ bellow $c_{2}$ there exist no states.
This inferior limit $c_{2}$ decreases with the increasing of $s$,
in the range $s\in (\pi /2, \pi )$, as prescribed in Sect.2 for type
\underline{(2)} and \underline{(2)} representations.

In Fig.6b we present the same potential for fixed $s=3.05$ 
($\simeq \pi $)and different
values for $m$.
We also remark here a very strong shape-dependence of the potential with
$m$.
For $m>2$ the potential is Morse-like and the depth of the valley increases
with $m$ like in the case of the undeformed limit [28].
For $m<2$ there are no more valleys and the potential is purely atractive.
In this case the spectrum is absolute continuous. This is again in agreement
with the analysis made in Sect.2. 

For the intermediate case $s=\pi /2$ we have $[2]=0$ and 
from eq.(54) and eq.(62) we obtain an exact harmonic oscillator
potential. If, in addition, $m\in Z$ and $F_{4}=0$ we obtain $V=const.$
This is also in perfect agreement with the consequences of the
analysis of the representations for this limiting case
$s \rightarrow \pi /2$ in Sect.2, ($e(2)$ dynamical
symmetry breaking). 
This situation
describes a 1-dimensional free particle and can be used similar with the
Euclidean connection [28], as the asymptotic limit for a scattering process.

By taking into account all these examples we can conclude that, the parameters
$s$ and $m$ provide the $su_{e^{is}}(2)$ model with a strong symmetry breaking effect:
the corresponding Schr\"odinger picture evoluates from periodic potential (Bloch-like)
made of negative  singularities into a periodic series of positive 
poles, or into
different bounded potentials allowing localised states (q-deformed
P\"oschl-Teller-like, Morse-like, harmonic oscillator like, etc.), or,
finaly, into constant potential (free asymptotic states).
This variety of the obtained shapes, when one modifies continuous $s$, and
discretly $m$, gives a wide area of possible 1-dimensional Schr\"odinger
pictures, all unified into the same dynamical symmetry, $su_{e^{is}}(2)$.
We note that similar results were obtained in the references [23,24],
concerning the deformation of an harmonic oscillator into a minus 
P\"oschl-Teller-like potential (similar with one of the cases
 presented in Fig.5b).
\par
\underline{General solutions}
\par

The general solutions of the system of eqs.(50,51) bring the most general
class of equivalent exact solvable potentials obtained through the
realisation (40,41). However, we mentione that there exist many other
possibilities of differential realisation of $su_{q}(2)$ (starting with
[5] up to a large variety of such realisations in the present literature).
In the above Subsection we have presented some limiting cases, when, due
to the decoupling of the two involved differential equations for the
function $f_{2}$ and for $R_{cm}$, four classes of potentials,
are obtained. These potentials are only combinations
of trigonometric ($cos(s)>0$), exponential ($s=\pi /2$) and
hyperbolic ($cos(s)<0$) functions. 

The general case, when the coefficient $f_{1}{\eta }^{2}[2m]$ in eq.(48)
is no more neglected, is much more interesting. Due to the coupling
of the differential equations for $f_{2}$ and the radial wavefunction
we need a special method of solving these problems.
First, $f_{1}$ is always given by eq.(49).
From eqs(48,50) we have
$$
R_{cm}^{'}=A(r)R_{cm},$$
\begin{eqnarray}
R_{cm}^{''}&=&A^{'}R_{cm}+AR_{cm}^{'},
\end{eqnarray}
with
\begin{eqnarray}
A&=&{{[2]} \over 2}[m]^{2}f_{1}-[2]\biggl (
{2 \over {{\eta }^{2}[2m]} }+{{[m]^{2}} \over {[2m]}}
\biggr )f_{2}-{4 \over {{\eta }^{2}[2m]}}{{f_{2}^{'}} \over {f_{1}}}.
\end{eqnarray}
By introducing $R_{cm}^{''}$ and $R_{cm}^{'}$ from eq.(64) in eq.(51),
$C(f_{1},f_{2},m,s)R_{cm}=cR_{cm}$, we obtain a non-homogenous,
non-linear differential equation with variable coefficients for $f_{2}$, in
the form
$$
a_{m}(r)f_{2}^{''}+
b_{m}(r)f_{2}^{'}+
c_{m}(r)(f_{2}^{'})^{2}+
d_{m}(r)f_{2}^{'}f_{2}+
$$
\begin{eqnarray}
e_{m}(r)f_{2}^{2}+
h_{m}(r)f_{2}+
g_{m}(r)&=&c
\end{eqnarray}
where the functions $a_{m}(r)$...$g_{m}(r)$ are obtained as combination
of the terms of $A$, $A^{2}$, $A^{'}$ and of those terms in the LHS of
eq.(51)
which do not contain the derivatives ${\partial }_{r}$.
The formal solution $f_{2}(r;m,s,c)$ of eq.(65), together with 
$f_{1}(r;m,s)$ from eq.(49), introduced in eq.(58) and with the help
of the substitution (61)
give a corresponding potential $V(r;m,s,c)$, and the eigenfunctions
\begin{eqnarray}
R_{cm}&=&R_{0}(c,m,s)exp \biggl (
\int \biggl (
{{[2]} \over 2}{\eta }^{2}[2m][m]^{2}f_{1}-[2]
f_{2}(2+{\eta }^{2}[m]^{2})-4{{f_{2}^{'}} \over {f_{1}}}
\biggr )dr
\biggr ).
\end{eqnarray}
The next step is to generate all the other states of the corresponding
unitary representations from which the first obtained solution, eq.(66),
belong, by acting with $J_{\pm }$ on it.
At present we do not know if this method allows the obtaining
of all the eigenstates, but
we can check its consistency at the algebra level.
Suppose we have calculated the solution $f_{2}$ of eq.(65)
and consequently we have obtained a first state given by eq.(66).
 We obtain a new
state in the form
\begin{eqnarray}
{\tilde R}_{cm}e^{i(m+1)}&=&J_{+cm}R_{cm}e^{im\phi },
\end{eqnarray}
where we denote with $J_{\pm cm}$, $C_{cm}$ the differential operators
eq.(40,41) in which we have introduced the solutions 
$f_{1,2}(r;m,s,c)$. We have
\begin{eqnarray}
C_{cm}{\tilde R}_{cm}e^{i(m+1)\phi }&=&c{\tilde R}_{cm}e^{i(m+1)\phi },
\end{eqnarray}
because $[C_{cm},J_{\pm cm}]=0$, and from eq.(2)
\begin{eqnarray}
[J_{+},J_{-}]_{cm}{\tilde R}_{cm}e^{i(m+1)\phi }&=&[2J_{z}]J_{+cm}
R_{cm}e^{im\phi }.
\end{eqnarray}
The eqs.(1,2) can be written in an equivalent form by 
using the example 1.5.3 given by Majid in [32]
$$q^{\pm 2J_{z}}J_{\pm }q^{\mp 2J_{z}}=q^{\pm 2}J_{\pm }.$$
We obtain from eqs.(69,70)
\begin{eqnarray}
[2J_{z}]J_{+cm}R_{cm}e^{im\phi }&=&[2(m+1)]{\tilde R}_{cm}
e^{i(m+1)\phi }.
\end{eqnarray}
From eqs.(69,71) we have then
$$J_{+cm}R_{cm}e^{im\phi }=R_{c,m+1}e^{i(m+1)\phi },
$$
which is in agreement with the action of $J_{+}$ in the representation $(c,m)$. 
The same thing happens with $J_{-}$.
In coordinates this reads
\begin{eqnarray}
{\Psi }_{c,m\pm 1}(r)&=&\biggl [
\pm R_{cm}^{'}(r)+\biggl (-{{f_{1}(r)} \over 2}[2m]+f_{2}\biggr )
R_{cm}(r)
\biggr ]e^{i(m \pm 1)\phi }.
\end{eqnarray}
\par
In conclusion of this section, the present 
method can generate different exact 1-dimensional
potentials and their corresponding eigenfunctions.
The problem of obtaining the full spectrum is strongly dependent
on the value of $s$, i.e. when $q$ is or is not a root of $1$.
 A complete investigation of such realisations is in course of
finalisation and it will be the subject of a following paper.
\vskip 1truecm
\underline{Remarks}
\vskip 0.5truecm
{\bf 1)} One interesting new point introduced here
in the differential realisation of
$su_{e^{is}}(2)$ 
is the resulted coupling between
$f_{2}$ and $R_{cm}$ in eq.(48). In the undeformed case  ($\eta =0$),
$f_{1,2}$ can be determined uniquely, and independent of $R_{cm}$.
Consequently one has a given $V(f_{1,2}(r),m)$ and can solve the
eigenproblem for this potential, i.e. to find the eigenstates according to
the unirreps of the undeformed algebra. More, since eq.(50), for $s=0$,
relates $f_{1}$ and $f_{2}$, if they are periodical functions, they should
have the same period. This case gives only periodic potentials with Bloch
structure, continuous spectrum (bands or full continuous) and delocalized
eigenstates.
The $q=e^{is}$ deformation allows, through the above mentioned coupling,
the occuring of new situations. Cases similar with the undeformed case
are encountered for $s=\pi {p \over l}$ ($p,l$ integers, $q$ a root of $1$),
when, for certain sets of $m$ the coupling term ${\eta }^{2}[2m]$
vanishes.
For $s\neq \pi {p \over l}$ ($q$ not a root of $1$), one can choose an
arbitrary function for $f_{2}$ and can solve $R_{cm}$ from eq.(50) only.
Then, eq.(51) is automatically fulfilled (since $R_{cm}e^{im\phi }$ belongs
to the space of the representations of $su_{e^{is}}(2)$ and $c$ results.
In such situations we can obtain, for a general potential
$V(f_{1,2}(r),m,s)$ only some of the solutions of the corresponding
Schr\"odinger equation.

If $s=\pi {p \over l}$ we obtain different solutions $R_{cm_{k}}$
for different $m_{k}$ such $[m_{k}]=[m_{k^{'}}],[2m_{k}]=[2m_{k^{'}}]$, ...
etc. If $s\neq \pi {p \over l}$ one can directly obtain only one solution
for each potential. The rest of the solutions for the same potentials are
obtained by the action of the generators $J_{\pm }$, as shown above.

We give an example: suppose $f_{1}$ is a periodical (of period
$\sqrt{cos(s)}$) and we look for the possibility 
of obtaining localised solutions $R_{cm}$.
Being localised, these solutions contain almost all the Fourier components
in the variable $r$ and, consequently, by solving eq.(48) with respect to
$f_{2}$, for such $R_{cm}$,
 we can obtain periodical solutions for $f_{2}(r)$ but in an incommensurate
period with respect to $\sqrt{cos(s)}$. As a consequence, the potential will
be a random potential and we have simulated, in this way, the Anderson
localization process [41].
On the contrary, if $R_{cm}$ is delocalized and consequently a periodic
function, by following the same methods $f_{2}$ will be also 
a periodic function and $f_{1,2}$, $V$
and $R_{cm}$ functions will have the same period, or integer multiples or
sub-multiples of the original one. In this case we have obtained the case of
commensurate potential, as one expects from the theory of solvable
potentials.
These two opposite aspects are controled through the $s$ parameter.
For example, when $q$ is a root of $1$ (special cyclic representations from
the point of view of q-groups [29-32]) there exist values for $m$ such that
$[2m]=0$. In this case the condition eqs.(48,50) ask for $f_{1}f_{2}=0$.
We can find a possible solution for this last condition, in the form
$$f_{1}(r)=\sum_{k\in I_{1}}C_{k}^{1}{\Phi }_{k}(r),
$$
\begin{eqnarray}
f_{2}(r)&=&\sum_{k\in I_{2}}C_{k}^{2}{\Phi }_{k}(r),
\end{eqnarray}
where $I_{1}\cap I_{2}=\{ 0\} $, $C_{k}^{1,2}$ are constants and
${\Phi }_{k}(r)=1$ for $k+\epsilon<r<(k+1)-\epsilon $, $k$ integer,
$0<\epsilon <1$,
and zero in the rest (wavelet functions). In this case $f_{1}$ and $f_{2}$
have disjoint support and the $f_{1}f_{2}=0$
condition is fulfilled.
For $\epsilon \rightarrow 1$, this example can provide exactly the
situation of periodic but incommensurable potentials [41].
\par
{\bf 2)} We want to note other two observations. First that the transition
from $su(2)$-like potentials (finite depth) towards $su(1,1)$-like potentials
(infinite wells) is performed through $su_{e^{is}}(2)$ in a continuous way,
by $s$. In this modality the q-deformed algebraic structure is a dynamical
symmetry for a larger class of potentials, including both limiting cases
(bounded and scattering states). Another interesting point arises when
looking at the turning point between 
these two different situations (change of sign of $cos(s)$). 
From the potential point of view this separations
in the two regions of deformation clasifies the potentials in
"trigonometric"-like and "hyperbolic"-like potentials, due to the solutions
of the functions $f_{1},f_{2}$ arrising in the differential realisation of
the generators. This condition is independent of the values of the quantum
numbers $c$ and $m$. By regarding at eq.(21) we remark that the same factor
$cos(s)$
occurs in the Casimir expression, i.e.: ${{\cos s} \over {\sin ^{2}s}}$,
which is the multiplicative factor in front of $\sin {2} (sJ_{z})$. 
So, $s$ can transform the Casimir element from a positive defined one
into an indefinite one.
For the case $\cos(s) < 0$ (which contains the $su(1,1)$ limit),
 we have, on one hand, hyperbolic-like
potentials (e.g. P\"oschl-Teller or deformation of it) which allow
both discrete and continuous spectra, infinite-dimensional. On the other
hand, in the same range for $s$, due to the fact that $c>{1 \over {\sin ^{2}
s}}$ , for any described series in Section 2, the "surfaces" described in
the end of Sect.2, defined by the condition $c=constant$, are open
surfaces (eq.(21) similar with the hyperboloidal surfaces connected with
SU(1,1). The reverse is also true: for
$cos(s) > 0$ (which contains the $su(2)$ limit) the resulting potentials are of trigonometric type
(depend
on trigonometric  functions) and have poles which divide the axis $r$
in infinite wells potential.
In this latter case the corresponding "topology" (Sect.2) is described
by an infinite reunion of compact spheres, similar with the $SU(2)$
topology.
This unexpected similarity between: the structure of the unitary 
representations, the
"shapes" associated with the Casimir operator eq.(19),
 and the behaviour of the
corresponding potentials, addres us the information and believe that
such approach of complex q-deformation of $su(2)$ is really a candidate
for the unification of the bounded and scattering systems in the same
dynamical symmetry.

\vskip 1.5truecm
\section{Conclusions}
\vskip 1truecm

In the present paper we have investigated some problems connected with the
complex deformation of the Lie algebra $su(2)$
with $q=e^{is}$ and through a general real deformation  of
 $su(2)$, as a Hopf algebra, too.

Using the usual trigonometric deformation [1-9,11-16,21-23]
for the q-algebra $su_{q}(2)$, we have constructed and clasified
some of
the coresponding unitary representations
 and we have shown the existence of the
possibility of an extension  of
the bound representations of $su(2)$ into three classes of
representations: continuous infinite-dimensional and discrete ones.
We note that the present $q=e^{is}$ deformations alow the "transition"
of some of the properties of the Lie algebra $su(2)$ (unirreps, Casimir 
operator)
towards those of the Lie algebras $su(1,1)$ and $e(2)$.
A pictorial illustration of the surfaces of constant value of the
q-deformed Casimir operator eigenvalues depending on $s$ and $c$
is described in a "formal" vector space of the
generators of $su(2)$. We note that even this example underline the
similarity between $su_{e^{is}}(2)$ and $su(1,1)$, for certain values of
$s$. 
Consequently there are chances for the application of such smooth
deformations in the AST theory, esspecially when one needs to connect
the bound and the scattering states in the same unifying (q-deformed
dynamical symmetry) picture.

We have proved the
existence of an exact analytical solution for 
a real generalised deformation functional 
$[2J_{z}]$, which
allows the deformation of the commutator relations of
$su(2)$ together with the introduction of a Hopf algebra structure.
This result (infinite dimensional spectra) shows that even for real
deformations one can modify the compact caracter of $su(2)$.
 
We have obtained a realisation of $su_{q}(2)$ in terms
of q-differential operators which could be related to an exact 
Schr\"{o}dinger
equation and a corresponding equivalent potential picture. 
We have solved the corresponding Schr\"odinger equation
in some limiting cases and we have shown some examples of deformed
potentials.
For some particular values of $s$ we have found
different shapes of potentials
of physical interest (solid state, nuclear, atomic and molecular interactions) 
like the examples given in Figs.4-6.
One can see from these examples that one can transform, through
the smooth variation of $s$, an infinite well potential into a
finite one, like was first suggested in [24].
\par
Consequently, due to the present approach, 
new potential shapes could get their algebraic
analogue in the corresponding nonlinear q-deformation.
The q-deformation of the bound states towards scattering states
can be applied to processes like scattering or decays 
in which the physical system completely changes its
symmetry. 
The deformed P\"oschl-Teller potential could
describe electrons in the valence band and the equivalent deformed
potential the electrons in the corresponding conduction band.
As an example one can use these new potentials also in the physics of
semiconductor junctions.  
 For example, the fact that, due to the coupling in eq.(48), the potential
depends on all the quantum numbers, $V(r;c,m,s)$, 
one could describe the dynamics of the "polaron"
in solid state structures or the deuteron states with one bound level
(self-localised states).
Direct applications in the theory of exact solvable 1-dimensional potentials
with random or incommensurate structure can also be found.
Since we have found some similarities between the rational/irrational
form of the deformation parameter $s/\pi $ with the
commensurate/incommensurate character of the resulting potentials,
we consider that exists a possibility to connect the $s$ parameter
with the order parameter of such potential 
models (which describes the random or
ordered statusof a lattice model.
In fact from eq.(40),
regarded as an eigenvalue problem,
 we can note a direct connection with driven-rotator
or other type of random oscillators, described by the eq.(2) in reference [41].
In this sens it will be interesting to further investigate the unitary
representations of $su_{e^{is}}(2)$ when $s/\pi $ is not rational
($q$ not a root of $1$) but belongs to a Liouville irrational class [41].
On the other hand, in the case when $q$ is a root of $1$ and the 
representations have
special distinguished properties and implication for the product of the
representations, we expect, via the above considerations, 
interesting posibilities of investigations of random/incommensurate
potentials for many-body problems. 
And further applications are
however, needed.

\vfill
\eject

\vfill
\eject
{\bf Figure captions}
\vskip 1truecm
\begin{enumerate}
\item[{\bf Fig.1}] Figures 1a-1b present pictorial views of some examples
of  
unitary representations of $su_{e^{is}}(2)$
for $s=1.013<\pi /2$. The horizontal axis contains the values of $m$
(integers and half-integers) and on the vertical axis contains 
the values of $c$, the
eigenvalue of the Casimir operator. We plot with horizontal lines the
limiting values for $c$, which separate the three main classes of unirreps,
i.e.:
$c_{0}={1 \over {\sin ^{2}s}},
c_{1}={1 \over {4\sin ^{2}(s/2)}}  c_{0}={1 \over {\cos ^{2}(s/2)}}$.
In the figures these classes of unirreps are labeled by the numbers (1)-(3) 
in the right part of each
plot, and the arrows indicate the ranges for $c$. 
The lower bounds of $c$, depending on $m$, (RHS of eq.(12)), 
are ploted with continuous lines denoted with 
{\bf 1} and {\bf 2} for the sign $\pm $, correspondingly.
functions.

\item[{\bf (a)}]       
Here we present the \underline{infinite dimensional unirreps}
(class \underline{(1)}), having their lower bound at $c=c_{0}$. 
The 14 full dots in the 
area $c>c_{0}$ are an example of such an unirrep for an arbitrary $c$.
Their range is marked with an arrow (1) in the right side of the figure.
In the same figure we plot
the \underline{mixed series,
(2)} in the range $c_{1}<c<c_{0}$ (see the arrow
(2) in the right side of the figure).
The full horizontal segments, having arrows at their ends,
represent the sequences of the
allowed/forbiden intervals for $m$, i.e. $J_{\delta }, J_{f}, 
J_{\Delta },...$, for a certain value of $c$, in the range $c\in (c_{0},
c_{1})$. We represent with full dots an example of an infinite-dimensional
unirrep belonging to class \underline{2a}, for $k=1$, and with stars a finite
dimensional unirrep for $N=4$, class \underline{2b}. 
Though the space between the values of $m$ is always 1,
in order to present these two situations in the same picture, we have used
two different values for $s$ (for the dots and the stars)but we 
have ploted them together, within the same 
scale.

\item[{\bf (a)}]
We present here the range
$c\in (c_{1},c_{2})$ for the \underline{discrete
unirreps}, class \underline{(3)} (see the arrow in the right side of the
figure). Again the full line with arrows shows the
allowed/forbiden intervals and the dots represent an example of such a
finit-dimensional unirrep for a certain $c$ in the range, and for $N=6$. 

\item[{\bf Fig.2}]
In this figure we present the dependence 
of the eigenvalues of $[J_{z}]$, $[2m]$, against the deformation parameter 
$s$.
One can see that starting from the undeformed case ($s=0$), after many
oscillations, in the final limit $s=\pi$, we obtain exactly the changing 
of the
sign, for the values of $m$, i.e. even values in this illustration,
 which exemplifies the similarity with 
the unirreps of $su(1,1)$. We note that for some special values of $s$, 
like $\pi /k=\pi /2, \pi /3, \pi /4$ etc, the level cross and we have
a pictorial view of  
what was calculated in Section 2, for the series \underline{2a}, case when 
there are exactely  $k$ distinct values for all $[2m]$.

\item[{\bf Fig.3}]
In this figure we present a cross section ($J_{y}=0$) 
in the formal vector space
$(J_{x},
J_{y},J_{z})$,  for the surfaces of constant value of the Casimir
operator ($J_{z}$-horizonat axis, $J_{x}$-vertical axis). 
This illustrates graphicaly the behaviour of the Casimir eigen
values for different valus of $s$. The 
three disjoint bigger circles represent the disconected
surfaces which appear in the case $\cos s > 0$ . 
The empty spaces
between these spheres result from the fact that the square root
which represent the expression of $J_{x}$ is complex. 
The central smaller sphere  represents the limiting case of $s=0$, i.e.
the spherical surface of the SU(2)
manifold.
For larger values of s, these spheres become "closer" one to each other and
when $\cos s < 0$  they join into one unique surface,
(the two pairs of continuous
ondulated open lines in figure, computed for two different values of 
$s\in (\pi /2, \pi )$), 
homotopical with the hyperboloidal shape associated with the SU(1,1)
manyfold. This $SU(1,1)$ hyperboloidal geometry is presented here
by the two extreme curves,
 each having  one single minimum (maximum)
respectively.

\item[{\bf Fig.4a and 4b}]
Potential shapes, ploted against $r$ in the case $cos(s)>0$, $m=1$.
having periodical structure. In Fig.4a one 
can see the difference between the limits
$\simeq 0$ and $s\simeq \pi /2$.
In Fig.4b we ploted for
fixed $s=0.25$ (rational, i.e. $q$ a root of $1$) and different $m$'s.
One can see the transition between negative and positive poles of
the potential, through the variation of $m$.

\item[{\bf Fig.5a and 5b}]
Potential shapes ploted against $r$ in the case $cos(s)<0$, for the first
set of solutions of eqs.(52,56).
 In Fig.5a we have fixed $m$ and variable $s$ in the range $\pi /2, \pi$.
One can see the transition from P\"oschl-Teller-like potentials
($s$ far from the limits) to (negative signed) harmonic oscillator-like potentials
($s\simeq \pi /2$, $s\simeq \pi $)
In Fig.5b for
fixed $s=3$ and different values for $m$, the dependence of the shape
with $m$ is presented. The difference between integer and half-integer
$m$ results in the symmetry of the shapes against $r=0$.

\item[{\bf Fig.6a and 6b}]
Potential shapes ploted against $r$, also in the case $s\in (\pi /2, \pi)$,
but for the second pair of solutions in eqs.(52,56).
In Fig 6a, for
fixed $m=1$ and variable $s\in (\pi /2, \pi )$ the potentials are similar
with
Morse potential. In Fig.6b for
fixed $s=3.05$ (closed to $\pi$) and different $m$ values.
the potentials transfor from Morse-like
($m>2$) into harmonic oscillator-like ($m<2$).

\end{enumerate}


\begin{thebibliography}{99}




\bibitem{1}
Drinfeld V G 1986  {\it Proc. Intern. Congress of
Mathematicians, Berkeley, CA, AMS, Providence, RI, p.798} 
( A. M. Gleason ed.);

Jimbo M 1985 {\it Lett. Math. Phys.} {\bf 10}, 63  

\bibitem{2}
Sklyanin E K 1982 {\it Funct. Anal. Appl.} {\bf 16} 262 

\bibitem{3}
Kulish P P and Reshetikhin Ny 1983 {\it J. Sov. Math.} {\bf 23}
2435

\bibitem{4}
Biedenharn L C 1990 {\it J. Phys. A: Math. Gen.} {\bf 22}
L873

\bibitem{5} 
Macfarlane A J 1990 {\it J. Phys. A: Math. Gen.} {\bf 22}
 4581
\bibitem{6}
Y J Ng 1990 {\it J. Phys. A: Math. Gen.} {\bf 23} 1023-1027

\bibitem{7}
Ui H and Aizawa N 1990 {\it Mod. Phys. Lett.} {\bf A 4} 237
\bibitem{8}
Kulish P P and Damaskinsky E V 1990 {\it J. Phys. A: Math. Gen.} {\bf 23}
L415
\bibitem{9} 
Maekawa T 1991 {\it J. Math. Phys.} {\bf 32} 2598;
T Masuda 1990 {\it et al.} {\it Lett. Math. Phys.} {\bf 19} 187.
\bibitem{10}
Gromov N A and Man'ko V I 1992 {\it J. Math. Phys.} {\bf 33}
1374
\bibitem{11} 
Pasquier V and Saleur H 1990 {\it Nucl. Phys. B} {\bf 330} 523
\bibitem{12} 
Bonatsos D Argyres E N, Drenska S B, Raychev P P,  
Roussev R P, and Smirnov Yu F 1990 {\it Phys. Lett. B} {\bf 251} 477
\bibitem{13}
 Cseh J, Gupta R K, Ludu A, Greiner W and Scheid W 1992 {\it J.
Phys. G: Nucl. Part. Phys.} {\bf 18} L73
\bibitem{14} Bonatsos D, Drenska S B, Raychev P P, Roussev R P, 
and
Smirnov Yu F 1992 {\it J. Phys. G: Nucl. Part} {\bf 17} L67
\bibitem{15} Bonatsos D, Raychev P P, Roussev R P, and Smirnov 
Yu F 1990
{\it Chem. Phys. Lett.} {\bf 175} 300 ;  
Chang Z and Yan H 1991 {\it Phys. Lett. A} {\bf 154} 254  
\bibitem{16}
Del Sol Mesa A, Loyola G, Moshinsky M and Velazquez V 1993
{\it J. Phys. A: Math. Gen.} {\bf 26} 1147
\bibitem{17}
Gupta R K and Ludu A 1993 {\it Phys. Rev.} {\bf C 48} 593
\bibitem{18}
Caracciolo R and Monteiro M R 1993 {\it Phys. Lett.} {\bf B 308} 58
\bibitem{19} 
Spiridonov V 1992 {\it Mod. Phys. Lett. A} {\bf 7} 1241;
S Skorik, V Spiridonov 1993 {\it Lett. Mat. Phys.} {\bf 28} 59-74. 
\bibitem{20} 
Spiridonov V 1992 {\it Phys. Rev. Lett.} {\bf 69} 398 
\bibitem{21} 
Fairlie D B and Nuyts J 1991 {\it J. Phys. A: Math. Gen.} {\bf 24} L1001
\bibitem{22}
Lerda A and Sciuto S 1993 {\it Nucl. Phys.} {\bf B 401} 613
\bibitem{23}
Bonatsos D, Daskaloyannis C and Kolokotronis P 1993 {\it J.Phys. A: Math.
gen.} {\bf 26} L871
\bibitem{24} 
Bonatsos D, Daskaloyannis C and Kokkotas K 1991 
{\it J. Phys. A: Math. Gen.} {\bf 24} L795; 
Bonatsos D, Daskaloyannis C and Kokkotas K 1992 {\it J. Math.
Phys.} {\bf 33} 2954
\bibitem{25} 
Polychronakos A P 1990 {\it Mod. Phys. Lett. A} {\bf 5} 2325; \\ 
Ro\v cek M 1991 {\it Phys. Lett. B} {\bf 255} 554 
\bibitem{26} 
Ludu A and Gupta R K 1993 {\it J. Math. Phys.} {\bf 34}  5367
\bibitem{27} 
Delbecq C and Quesne C 1993 {\it J. Phys. A: Math. Gen.} {\bf 26} L127; 
1993 {\it Phys. Lett. B} {\bf 300} 227; 
{\it Mod. Phys. Lett. A} 
{\bf 8} 961; 
Bonatsos D, Daskaloyannis C, Kolokotronis P, Ludu A and Quesne C 1995
{\it Lett. Math. Phys.}, submitted
\bibitem{28}
Alhassid Y, G\"ursey F and Iachello F 1983 {\it Annals of Phys.} {\bf 148}
346; \\        
Alhassid Y, G\"ursey F and Iachello F 1986 {\it Annals of Phys.} {\bf 167} 
181; \\    

\bibitem{29}
L. C. Biedernharn 1990 {\it An Overview of Quantum Groups}
(The XVIII$^{th}$ International Conference on Group Theoretical Methods
in Physics, Moscow, June 4-9). 
   
\bibitem{30}
 Z. Chang 1995 {\it Phys. Rep.} {\bf 262} 137-125;
de Concini C and Procesi C 1993 in {\it D-modules, Representation Theory
and Quantum Groups} (Springer-verlag, Zampieri G and D'Agnolo A, eds.).

\bibitem{31}
A M Gavrilik, A U Klimyk 1994 {\it J. Math. Phys.} {\bf 35} 3670-3686.

\bibitem{32}
Majid S 1990 {\it Int. J. Mod. Phys. A} {\bf 5} 1; 
Abe E 1980 {\it Hopf Algebras} (Cambridge Tracts in
Mathematics, Cambridge)
\bibitem{33} 
Dobrev V K April 1991 {\it Introduction to Quantum Groups}( preprint
Institute for Theoretical Physics, University of Gottingen); \\
Alexeev A, Gluschenkov D and Lyakhovskaya A June 1992 {\it Preprint
St.Petersburg State University, PAR-LPTHE 92-94};
1992 {\it Lectures on Quantum Groups} in {\it Quantum Groups, Integrable
Statistical Models and Knot Theory} (World Scientific, Ge L M and de Vega
H S , eds.)

\bibitem{34}
Curtright T L and Zachos C K 1990 {\it Phys. Lett} {\bf B 243} 237
\bibitem{35}
Wybourne B G 1974 {\it Classical groups for Physicists} (Wiley, New York)
\bibitem{36}
Celeghini E, Giachetti E, Sorace E and Tarlini M 1990 {\it J. Math. Phys.}
{\bf 31} 2548
\bibitem{37}
Polychronakos A P 1990 {\it Mod.Phys.Lett A} {\bf 5} 2325
\bibitem{38} 
Daskaloyannis C 1992 {\it J. Phys. A: Math. Gen} {\bf 25}
2261
\bibitem{39} 
Carlson B C 1977 {\it Special Functions of Applied Mathematics}
(Academic Press New York)
\bibitem{40} 
Zielke A 1993 {\it Doctor Thesis}, Giessen University, unpublished;
J Wu 1991 {\it Doctor Thesis}, Yale University, unpublished.
\bibitem{41}
Grempel D R, Fishman S and Prange R E 1982 {\it Phys. Rev. Lett}
{\bf 49} 833; 1983 {\it Phys. Rev. B} {\bf 28} 7370; 1984 {\it Phys. Rev. B}
{\bf 29} 6500 and Simon B 1985 {\it Ann. Phys.} {\bf 159} 157.

\end{thebibliography}
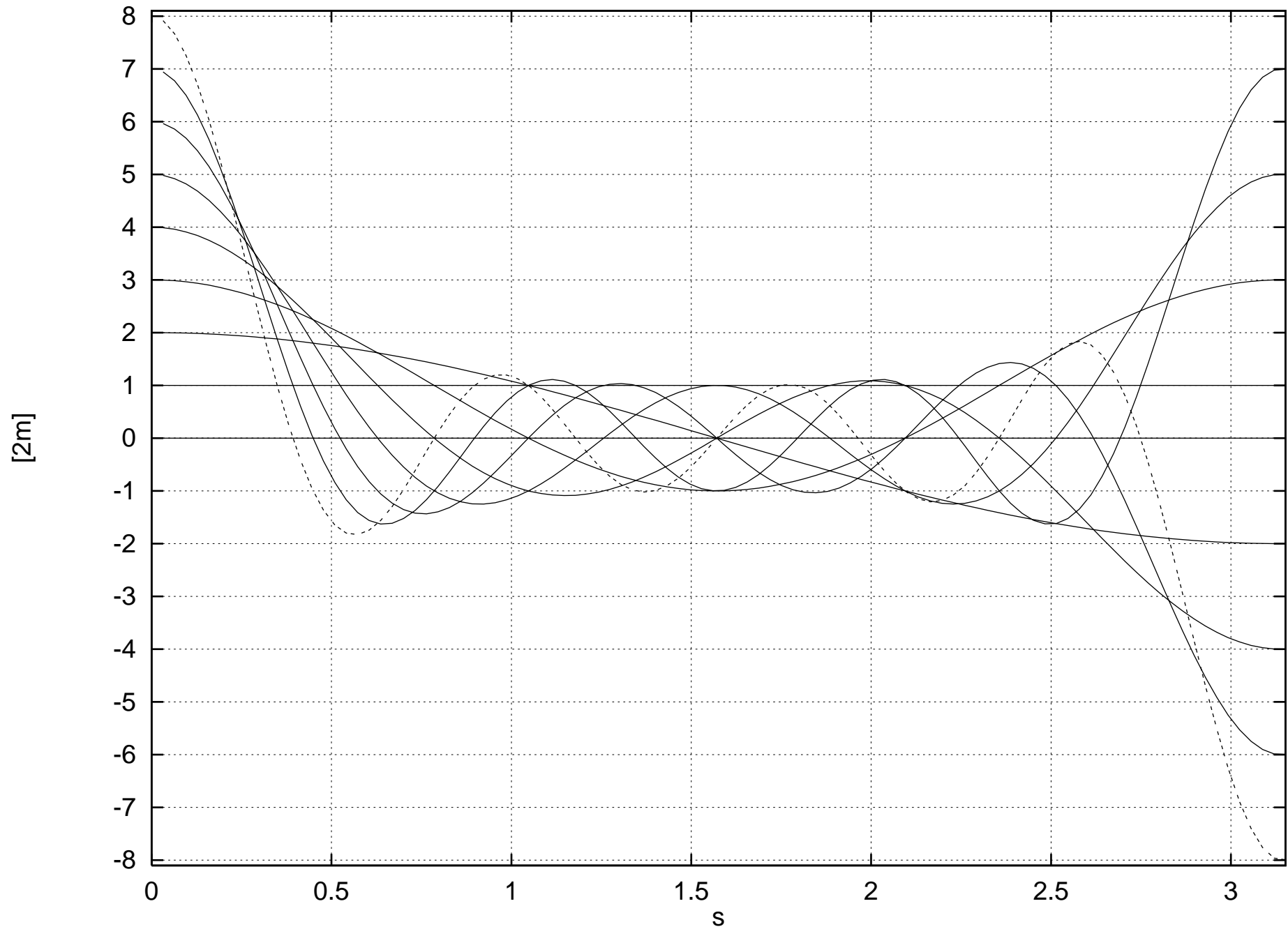
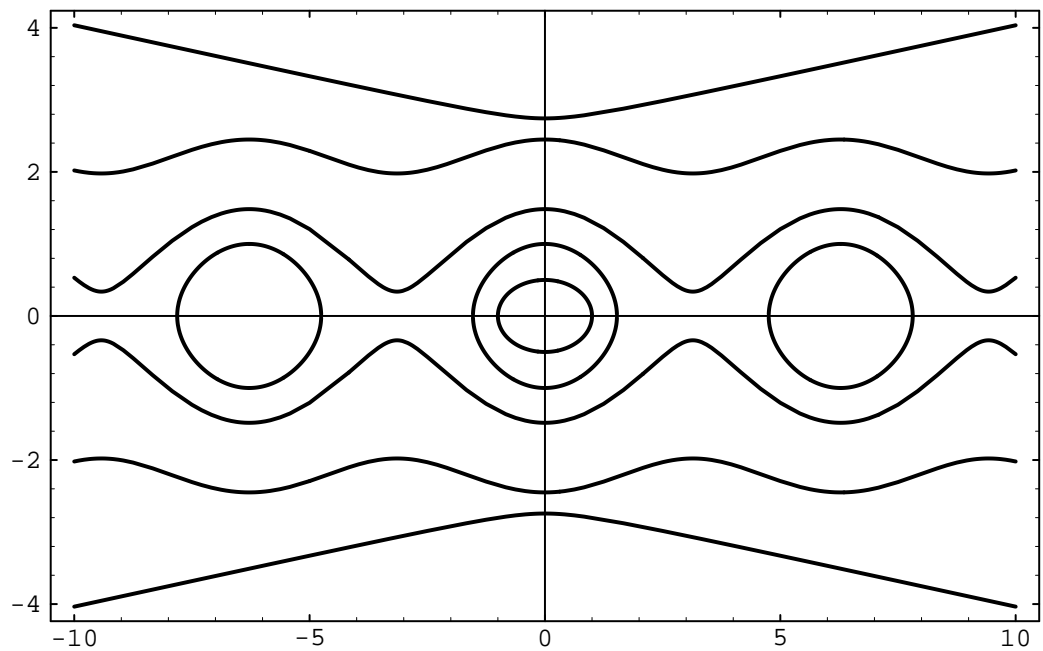
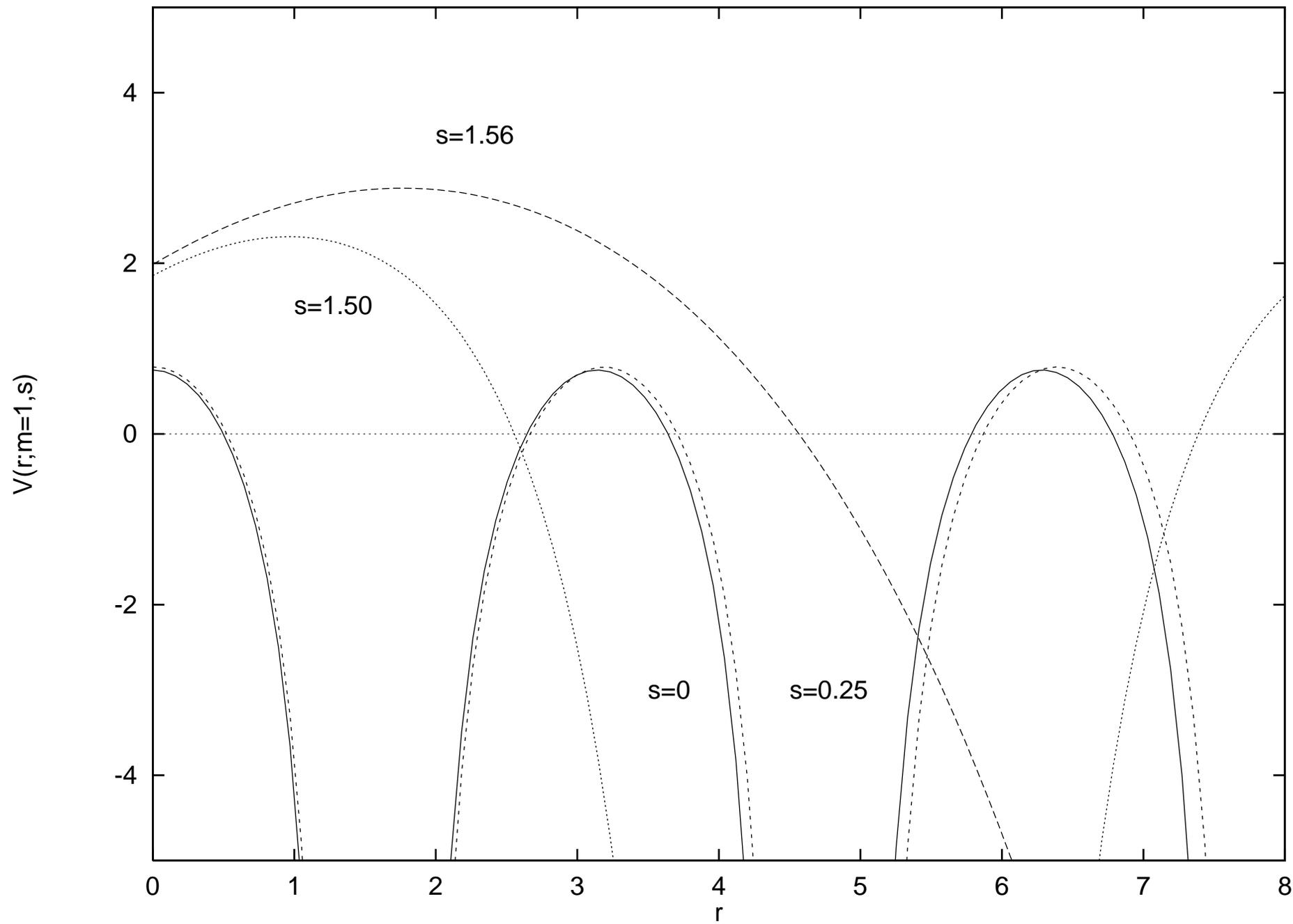
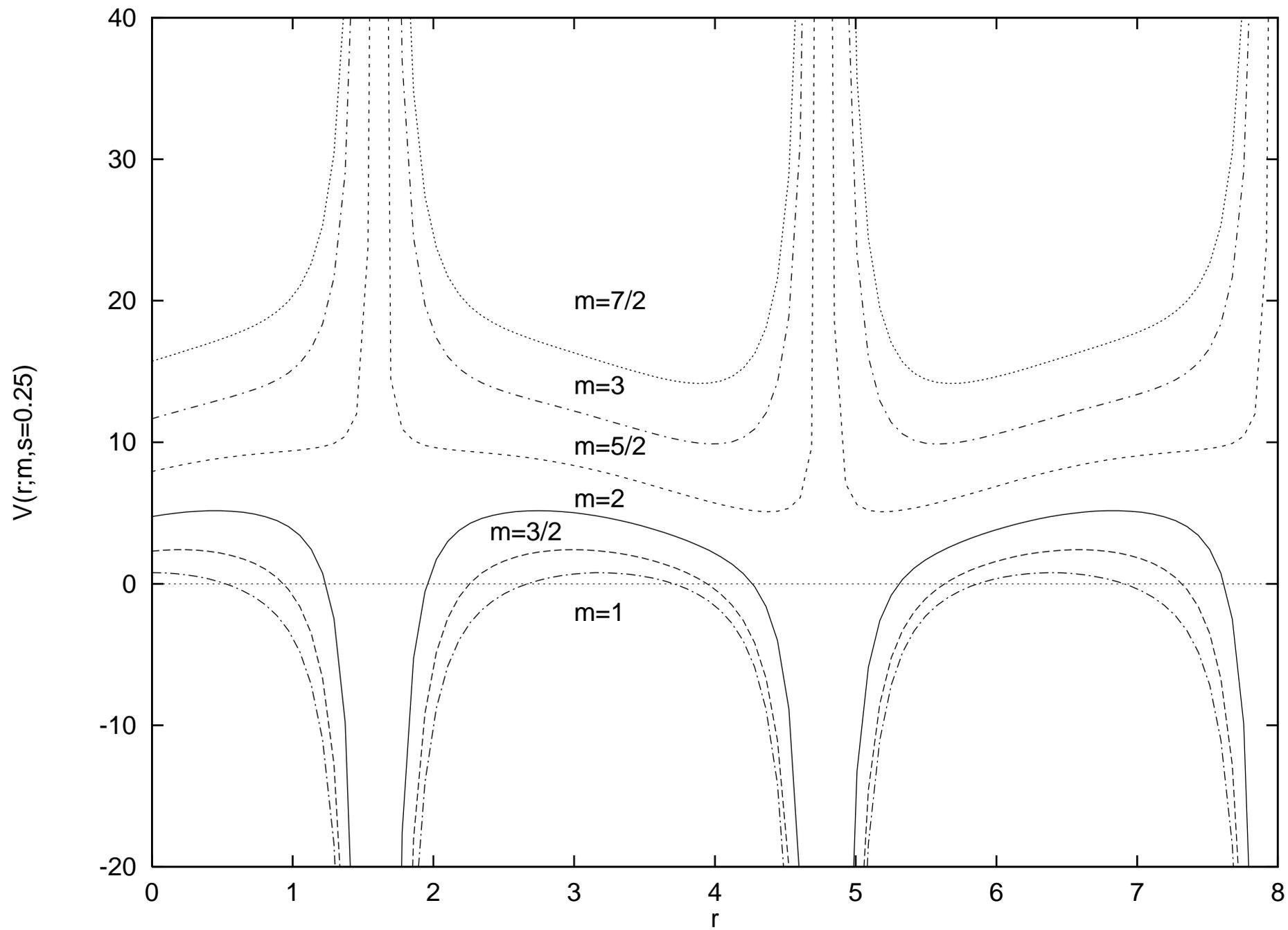
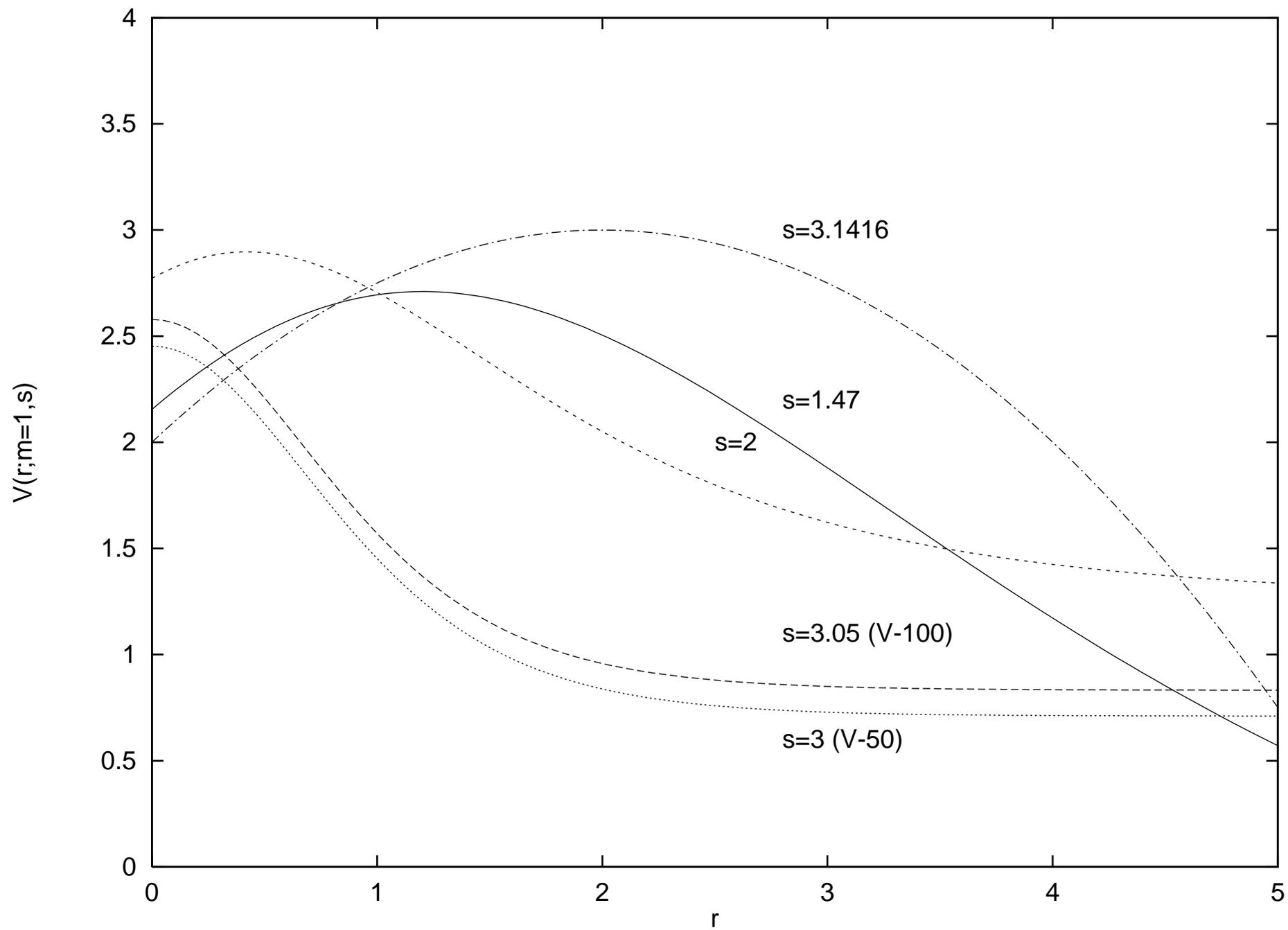
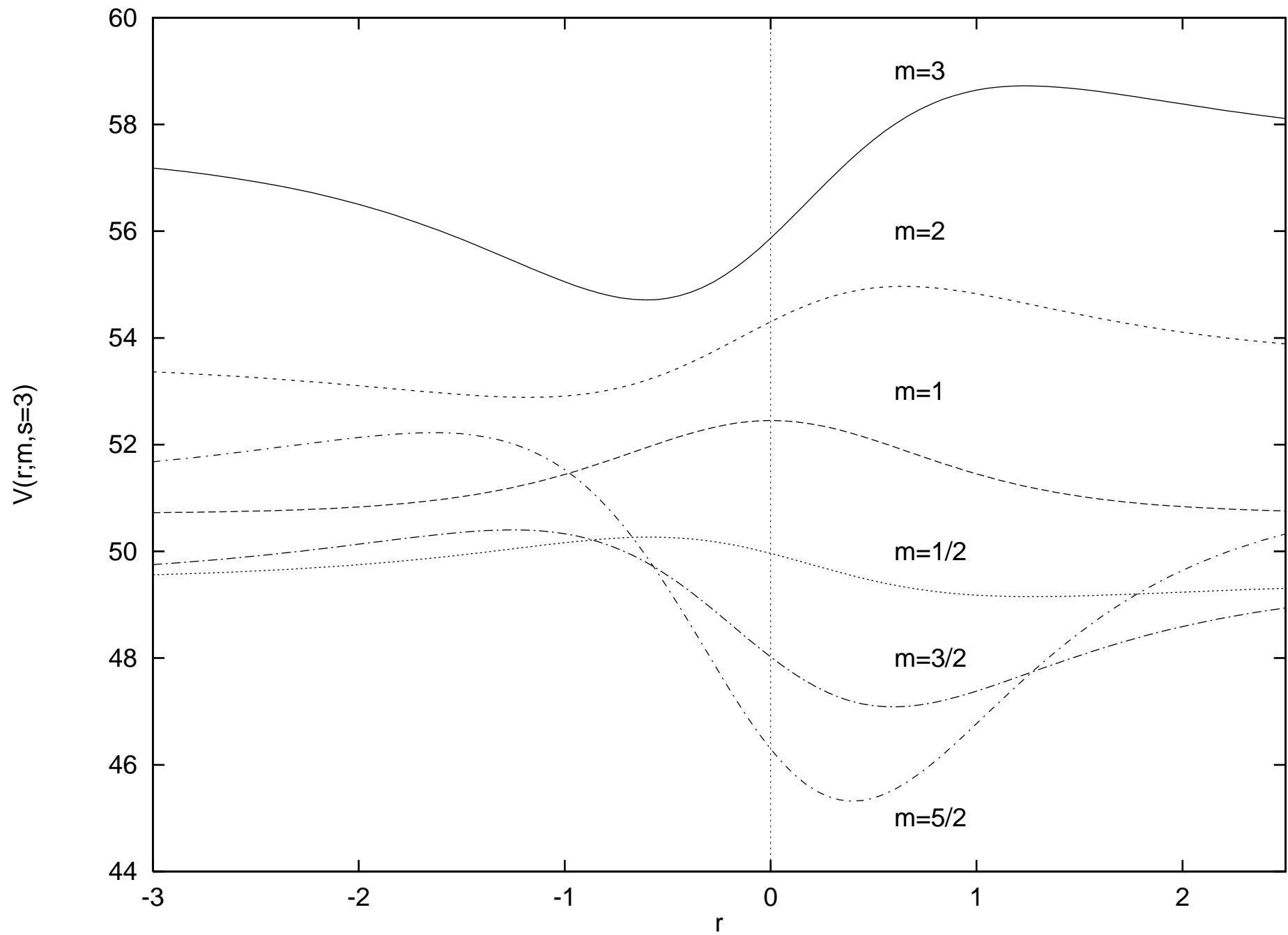
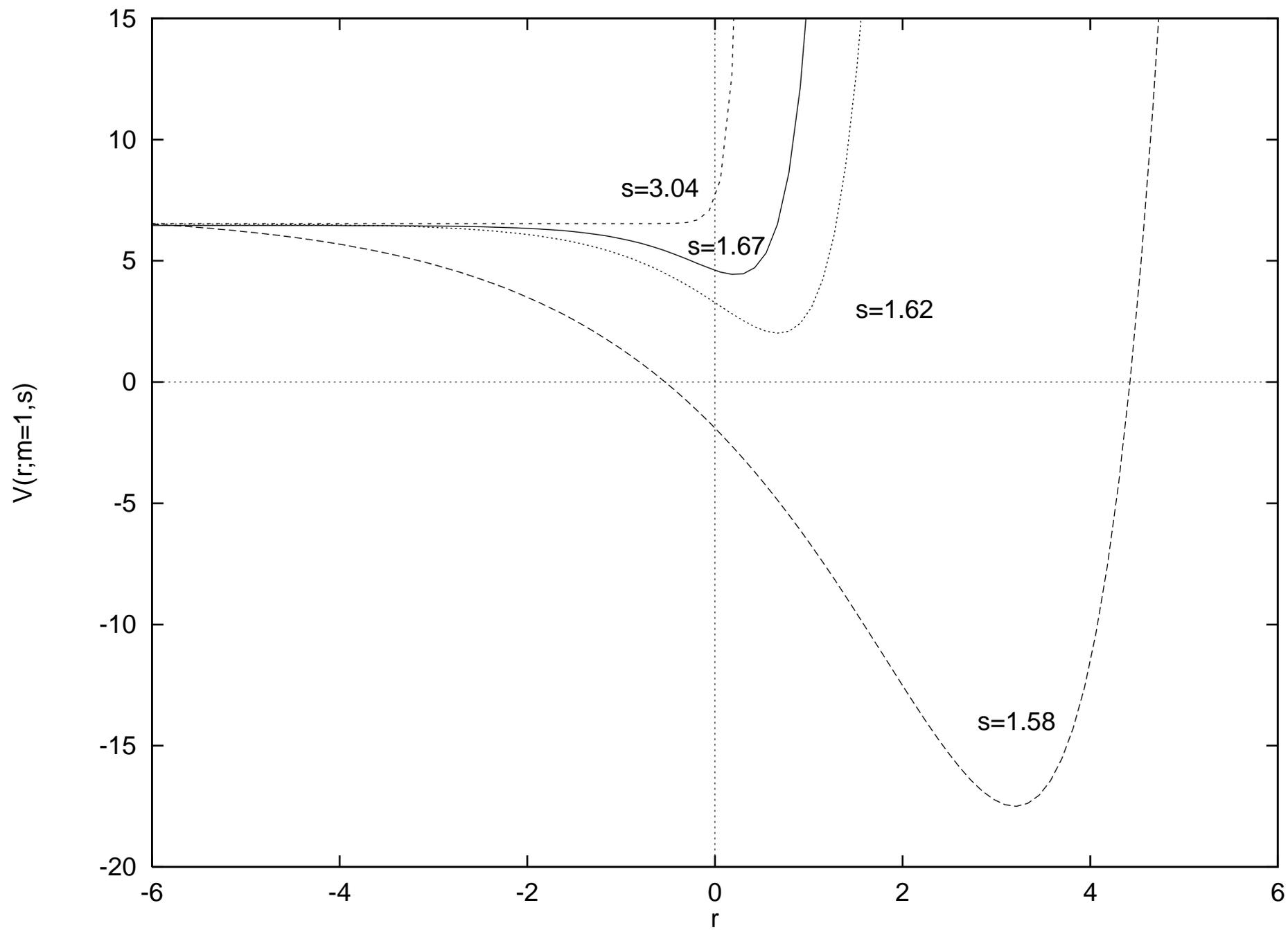
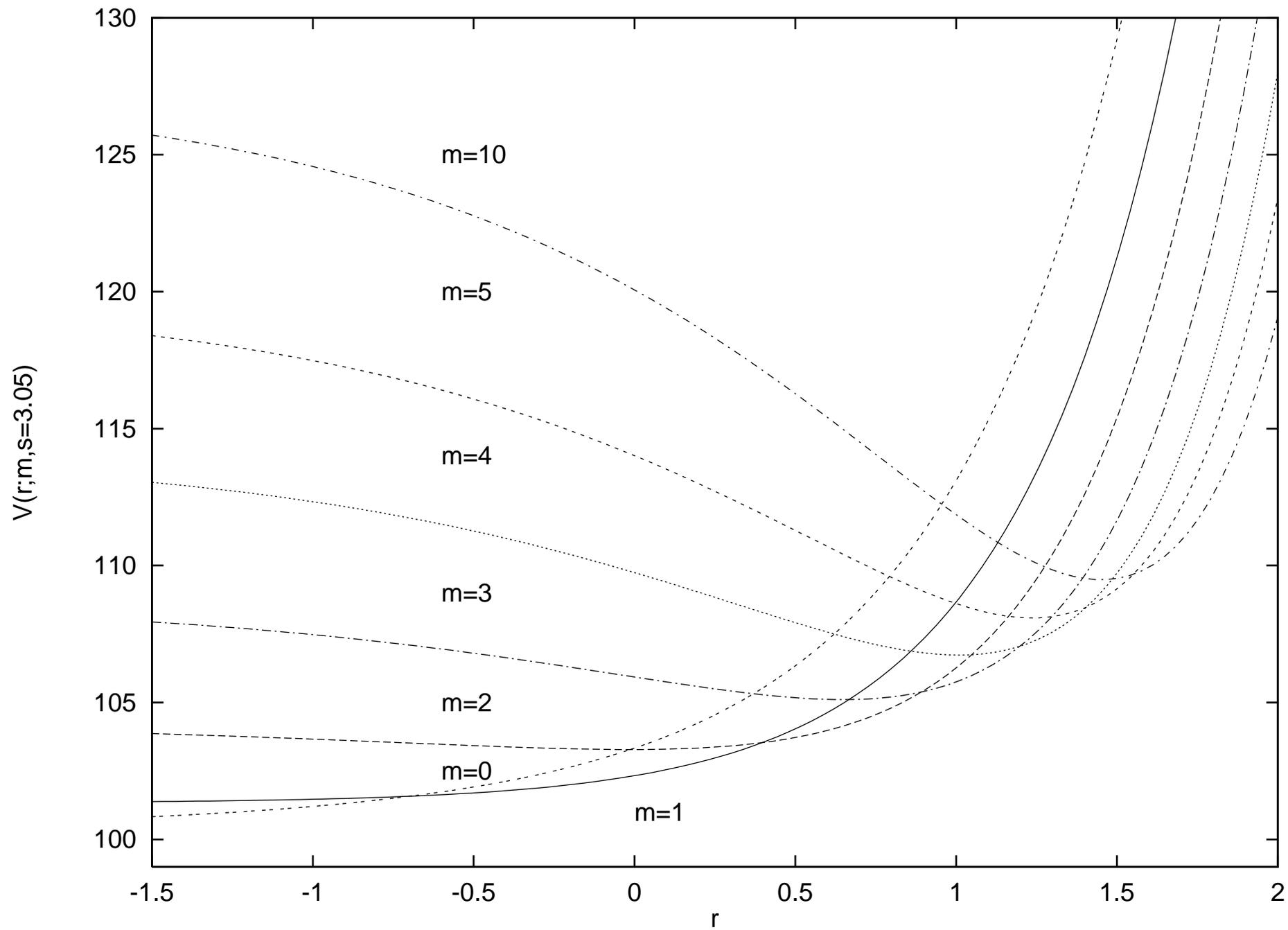
\end{document}